\documentclass[onefignum,onetabnum]{siamonline250211}
\usepackage{filecontents}
\usepackage{pdfpages}
\usepackage{amssymb}

\usepackage{booktabs}

\usepackage{lipsum}
\usepackage{amsfonts}
\usepackage{graphicx}
\usepackage{epstopdf}
\usepackage{mathtools}
\ifpdf
  \DeclareGraphicsExtensions{.eps,.pdf,.png,.jpg}
\else
  \DeclareGraphicsExtensions{.eps}
\fi
\usepackage{algorithm}
\usepackage[noEnd,indLines,commentColor=gray]{algpseudocodex}
\usepackage{enumitem}
\setlist[enumerate]{leftmargin=.5in}
\setlist[itemize]{leftmargin=.5in}

\newsiamremark{remark}{Remark}
\newsiamremark{hypothesis}{Hypothesis}
\crefname{hypothesis}{Hypothesis}{Hypotheses}
\newsiamthm{claim}{Claim}
\newsiamremark{fact}{Fact}
\crefname{fact}{Fact}{Facts}
\newsiamremark{assumption}{Assumption}
\crefname{assumption}{assumption}{assumptions}   
\Crefname{assumption}{Assumption}{Assumptions}   

\headers{Inertial Langevin Algorithm}{A. Falk, A. Habring, C. Griesbacher, and T. Pock}

\title{Accelerated Sampling with the Inertial Langevin Algorithm\thanks{Initial version submitted to the editors on December 9, 2025. Revised version submitted on August 31, 2026.
\funding{This research was funded in whole or in part by the Austrian Science Fund (FWF) (\href{https://www.fwf.ac.at/en/research-radar/10.55776/COE12}{10.55776/COE12}, \href{https://www.fwf.ac.at/en/research-radar/10.55776/F100800}{10.55776/F100800}) and the Lead project \href{https://www.tugraz.at/en/news/article/tu-graz-buendelt-ihre-kraefte-in-biotechnologie-und-kuenstlicher-intelligenz}{DigiBioTech} of Graz University of Technology.}}}

\author{%
    Alexander Falk\thanks{%
        Institute of Visual Computing, Graz University of Technology (\email{falk@tugraz.at}, \email{andreas.habring@tugraz.at}, \email{griesbacher@tugraz.at}, \email{thomas.pock@tugraz.at}).%
    }
    \and Andreas Habring\footnotemark[2]
    \and Christoph Griesbacher\footnotemark[2]
    \and Thomas Pock\footnotemark[2]%
}

\usepackage{amsopn}
\DeclareMathOperator{\diag}{diag}

\usepackage{siunitx}
\usepackage{printlen}

\usepackage{bbm}

\DeclareMathOperator{\id}{I}
\newcommand{\expo}[1]{\mathrm{e}^{#1}}

\DeclareMathOperator{\TV}{TV}

\DeclareMathOperator{\trace}{tr}
\DeclareMathOperator*{\argmin}{argmin}

\DeclareMathOperator{\cov}{cov}

\DeclareMathOperator{\law}{law}

\newcommand{\doublehookrightarrow}%
{\DOTSB\lhook\joinrel\relbar\!\!\!\!\lhook\joinrel\rightarrow}

\newcommand{\longrightharpoonup}%
{\relbar\joinrel\rightharpoonup}

\newcommand{\conditionalcomma}[1]{\ifx#1\empty\else,\fi}

\newcommand{\Wc}{\mathcal{W}}
\newcommand{\Pc}{\mathcal{P}}

\newcommand{\Vc}{\mathcal{V}}

\newcommand{\1}{\mathbbm{1}}
\newcommand{\R}{\mathbb{R}}
\newcommand{\E}{\mathbb{E}}
\renewcommand{\P}{\mathbb{P}}
\newcommand{\N}{\mathbb{N}}

\newcommand{\poly}{\mathrm{I\mspace{-2.5mu}P}}

\newcommand{\grad}{\nabla}

\newcommand{\abs}[1]{{|{#1}|}}

\newcommand{\dd}{\mathrm{d}}

\newcommand{\iid}{i.i.d.}
\renewcommand{\epsilon}{\varepsilon}

\newcommand{\Bc}{\mathcal{B}}

\newlength{\formulaindentwidth}

\newcommand{\Nc}{\mathcal{N}}

\newcommand{\Oc}{\mathcal{O}}

\newcommand{\step}{{\Delta t}}
\newcommand{\stepd}{{\tau}}
\newcommand{\disca}{{\rho}}

\newcommand{\friction}{\varepsilon}

\newcommand{\wiener}{{\dd W_t}}

\renewcommand{\phi}{\varphi}

\newcommand{\inner}[2]{\left\langle#1,#2\right\rangle}

\newcommand{\grc}{\alpha}

\renewcommand{\poly}{\mathrm{poly}}

\usepackage{glossaries-extra}
\newabbreviation{ula}{ULA}{unadjusted Langevin algorithm}
\newabbreviation{ald}{ALD}{annealed Langevin dynamics}
\newabbreviation{myula}{MYULA}{Moreau-Yosida regularized unadjusted Langevin algorithm}
\newabbreviation{daz}{DAZ}{diffusion at absolute zero}
\newabbreviation{apgd}{APGD}{accelerated proximal gradient descent}
\newabbreviation{dsm}{DSM}{denoising score matching}
\newabbreviation{sde}{SDE}{stochastic differential equation}
\newabbreviation{map}{MAP}{maximum a-posteriori}
\newabbreviation{mcmc}{MCMC}{Markov chain Monte Carlo}
\newabbreviation{mc}{MC}{Markov chain}
\newabbreviation{mmse}{MMSE}{minimum mean-squared-error}
\newabbreviation{tdv}{TDV}{total deep variation}
\newabbreviation{mri}{MRI}{magnetic resonance imaging}
\newabbreviation{iid}{i.i.d.}{independent and identically distributed}
\newabbreviation{lsc}{l.sc.}{lower semicontinuous}
\newabbreviation{tv}{TV}{total variation}
\newabbreviation{bp}{BP}{belief propagation}
\newabbreviation{gmm}{GMM}{Gaussian mixture model}
\newabbreviation{kde}{KDE}{kernel density estimation}
\newabbreviation{kld}{KLD}{Kullback-Leibler divergence}
\newabbreviation{tvd}{TVD}{total variation distance}
\newabbreviation{sor}{SOR}{successive over-relaxation}
\newabbreviation{rv}{RV}{random variable}
\newabbreviation{ila}{ILA}{inertial Langevin algorithm}
\newabbreviation{em}{EM}{Euler-Maruyama}
\newabbreviation{nila}{NILA}{Nesterov's inertial Langevin algorithm}
\newabbreviation{skrock}{SK-ROCK}{SK-ROCK}
\newabbreviation{ou}{OU}{Ornstein-Uhlenbeck}
\newabbreviation{md}{MD}{molecular dynamics}
\newabbreviation{gd}{GD}{gradient descent}
\newabbreviation{hmc}{HMC}{Hamiltonian Monte Carlo}
\newabbreviation{ses}{SES}{stochastic exponential Euler scheme}
\newabbreviation{emd}{EMD}{earth mover's distance}
\newabbreviation{acf}{ACF}{autocorrelation function}
\newabbreviation{ot}{OT}{optimal transport}
\newabbreviation{rof}{ROF}{Rudin-Osher-Fatemi}
\newabbreviation{glm}{GLM}{Gaussian latent machine}
\newabbreviation{lm}{LM}{Leimkuhler-Matthews}
\newabbreviation{ode}{ODE}{ordinary differential equation}
\newabbreviation{ebm}{EBM}{energy-based model}
\newabbreviation{pes}{PES}{potential energy surface}
\newabbreviation{egnn}{EGNN}{equivariant graph neural network}
\glsdisablehyper

\newcommand{\ie}{\textit{i.e.}}
\newcommand{\eg}{\textit{e.g.}}
\newcommand{\cf}{\textit{cf.}}

\DeclareSIUnit\angstrom{\text {Å}}

\ifpdf
\hypersetup{
  pdftitle={Inertial Langevin Algorithm},
  pdfauthor={A. Falk, A. Habring, C. Griesbacher, and T. Pock}
}
\fi

\Crefname{SM:eq}{section}{sections}

\usepackage{thmtools}
\usepackage{hyperref}
\usepackage{cleveref}

\begin{document}

\maketitle


\begin{abstract}
We consider the \emph{inertial Langevin algorithm} (ILA), a simple momentum-based method for sampling from Gibbs distributions of the form $\pi(x) \propto \exp(-U(x))$.
\glsunset{ila}\Gls{ila} augments the unadjusted Langevin algorithm with an inertia term and a matching noise rescaling, yielding a sampling analogue of Polyak's heavy-ball method from optimization.
This modification provably accelerates convergence.
For Gaussian targets, ILA attains Wasserstein-2 accuracy $\delta$ in $\tilde{\mathcal O}(\sqrt{\kappa}/\sqrt\delta)$ iterations, where $\kappa = L/m$ denotes the condition number and $\tilde{\mathcal O}$ hides logarithmic factors, thus reaching the ballistic $\sqrt{\kappa}$ complexity known from accelerated optimization and proven to be a lower bound for sampling from Gaussians~\cite{chewi2024query}. Moreover, we show that the method's bias can be fully removed and, thus, complexity reduced to $\tilde \Oc(\sqrt{\kappa})$ by considering the sequence consisting of the average of ever two consecutive \gls{ila} iterates, which is shown to be equivalent to the well-known BAOAB splitting scheme in the linear setting.
Beyond Gaussians, for $L$-smooth potentials with sufficient growth at infinity, we prove geometric ergodicity of ILA and of its continuous-time analogue, the underdamped Langevin dynamics, together with convergence of the discretization bias to zero as the step size decreases.
These guarantees hold under considerably simpler parameter restrictions than previously available in the literature, in particular, not imposing a lower bound on the friction, which is conjectured to be crucial for acceleration~\cite{altschuler2026shifted}.
We underpin our theoretical findings with numerical experiments covering ill-conditioned Gaussian distributions, total variation image denoising, and the challenging task of maximum likelihood learning of an energy-based model for molecular structure generation.
\end{abstract}



\section{Introduction}
\begin{figure}
    \centering
   \includegraphics[width=1.0\linewidth]{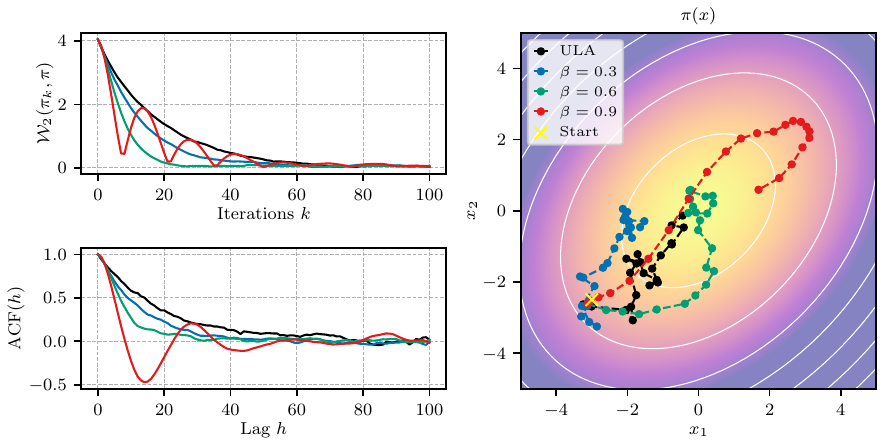}
   \vspace{-0.75cm}
    \caption{Sweep over momentum parameter $\beta$ on a simple 2D Gaussian distribution.
    \emph{Right}: It can be seen that for increasing $\beta$ the process performs larger steps, leading to better exploration of the sample space.
    \emph{Left top:} This behavior also manifests on a distribution level, where larger momenta lead to faster convergence in Wasserstein-2 distance.
    It is also observable that setting $\beta$ too large leads to oscillatory behavior characteristic of momentum-based optimization methods.
    \emph{Left bottom:} With increasing momentum, successive samples within a chain become more widely spaced.
    Consequently, the \gls{acf} decays faster, yielding larger effective sample sizes from a single chain.
    }
    \label{fig:momentum_sweep}
\end{figure}
In this work, we address the problem of sampling from Gibbs probability distributions on $\R^d, d\geq1$ with Lebesgue density\footnote{In a slight abuse of notation, depending on the context, $\pi$ may denote either the density or the measure itself.}
\begin{equation}\label{eq:gt-density}
    \pi (x) = \frac{\exp(-U(x))}{\int_{\R^d} \exp(-U(\xi)) \, \dd \xi},
\end{equation}
with a potential $U \colon \R^d \to \R$.
This is a task frequently occurring, \eg, in Bayesian inverse problems~\cite{luo2023bayesian}, mathematical imaging~\cite{durmus2018efficient,pereyra2016proximal}, machine learning~\cite{zach2023stable,zach2021computed}, and uncertainty quantification~\cite{narnhofer2022posterior}.
We consider the~\glsfirst{ila}, a gradient-based \gls{mcmc} sampler that, starting from $X_1 = X_0$, iterates
\begin{equation}\label{eq:ila}
    X^{k+1} = X^k - \stepd \grad U(X^k) + \beta (X^{k} - X^{k-1}) + \sqrt{2 \stepd (1 - \beta)} N^k,\quad k=1,2,\dots,
    \tag{ILA}
\end{equation}
where $\tau > 0$ denotes the step size, $\beta \in [0,1)$ is a \emph{momentum} parameter and $(N^k)_k$ is a sequence of \iid~standard Gaussian vectors in $\R^d$.
\eqref{eq:ila} is a stochastic counterpart of Polyak's heavy-ball algorithm~\cite{polyak1964some} for minimizing $U$, which reads as
\begin{equation}
    x^{k+1} = x^k - \tau \nabla U(x^k) + \beta (x^k-x^{k-1}), \quad k=1,2,\dots.
\end{equation}
Similar to heavy-ball being an adaptation of mere gradient descent,~\cref{eq:ila} is an adaptation of the \gls{ula}~\cite{robert1999monte,durmus2017nonasymptotic,durmus2019high,dalalyanTheoreticalGuaranteesApproximate2017a,dalalyan2017further},
\begin{equation}
    X^{k+1} = X^k - \stepd \grad U(X^k) + \sqrt{2 \stepd}N^k, \quad k=1,2,\dots,
    \tag{ULA}
\end{equation}
which is evidently recovered from \cref{eq:ila} by setting $\beta=0$.
\Gls{ula} is the \gls{em} discretization of the overdamped Langevin \gls{sde}
\begin{equation}\label{eq:OD}
    \dd X_t = -\nabla U(X_t) \, \dd t + \sqrt{2} \, \dd W_t,\tag{OLD}
\end{equation}
where $(W_t)_t$ is a standard $d$-dimensional Brownian motion.
\Cref{fig:momentum_sweep} illustrates the effect of momentum on the sampled trajectories.

It is well-known in optimization, that, for quadratic potentials with condition number $\kappa$\footnote{The condition number $\kappa$ is the quotient $L/m$ when $L$ and $m$ are the largest and smallest eigenvalues of $\nabla^2 U$, respectively.}, the iteration complexity of the heavy ball method improves to $\Oc(\sqrt{\kappa})$, compared with $\Oc(\kappa)$ for ordinary gradient descent.
The heavy-ball method is itself a discretization of the gradient flow with friction
\begin{equation}\label{eq:ode}
    \ddot x = -\nabla U(x) - \friction \dot x,
\end{equation}
whose stochastic analog is the \emph{underdamped} (or \emph{kinetic}) Langevin process~\cite{maThereAnalogNesterov2021}, defined as
\begin{equation}\label{eq:uld_intro}
    \begin{cases}
        \dd X_{t} = V_t \, \dd t,\\
        \dd V_t = [-\gamma V_t - \nabla U(X_t)] \,\dd t + \sqrt{2\gamma} \, \wiener.\tag{ULD}
    \end{cases}
\end{equation}
Under appropriate conditions on $U$,~\eqref{eq:uld_intro} is ergodic in a variety of metrics with unique stationary distribution $p(x,v)\propto \exp(-U(x) - \frac{\|v\|^2}{2})$~\cite{cheng2018underdamped,cao2023explicit,dalalyan2020sampling,durmusUniformMinorizationCondition2025,gouraud_hmc_2025,leimkuhler2024kinetic,maThereAnalogNesterov2021,pavliotis2014stochastic,schuh2024convergence}.
Formally, depending on the choice of $\beta$, \gls{ila} can be understood either as a discretization of~\cref{eq:OD} or of~\cref{eq:uld_intro}:
For fixed $\beta$,~\eqref{eq:ila} yields
\begin{equation}\label{eq:overdamped_limit_ila}
    \frac{X^{k+1}-X^k}{\tau} - \beta\frac{X^k-X^{k-1}}{\tau} = -\nabla U(X^k) + \frac{\sqrt{2(1-\beta)}}{\sqrt{\tau}}N^k,
\end{equation}
which for $\tau \to 0$ informally converges to
\begin{equation}\label{eq:ila_to_ula}
    \dot X = -(1-\beta)^{-1}\nabla U(X) + \sqrt{2(1-\beta)^{-1}}\dot W,    
\end{equation}
\ie, overdamped Langevin dynamics preconditioned with $(1-\beta)$. 
Therefore, larger $\beta$ accelerates the continuous-time dynamics.
This derivation serves also as a motivation for the noise scaling in~\eqref{eq:ila} as it ensures that~\cref{eq:ila_to_ula} preserves $\pi$ as its stationary distribution.
On the other hand-side, the Gaussian analysis below suggests the $\tau$-dependent choice $\beta = (1-\friction\sqrt{\tau})^2 + \Oc(\tau)$ for some $\friction>0$.
With the auxiliary variable $V^k \coloneqq (X^k-X^{k-1})/\sqrt{\tau}$ we have
\begin{equation}\label{eq:ILA_is_underdamped}
    \begin{cases}
        V^{k+1} =& V^k - 2\friction\sqrt{\tau}V^k + \Oc(\tau) V^k -\sqrt{\tau} \nabla U(X^k) + \sqrt{2\friction\sqrt{\tau}}N^k\\
        X^{k+1} =& X^k + \sqrt{\tau} V^k
    \end{cases}
\end{equation}
which is easily seen to be a discretization of~\eqref{eq:uld_intro} upon identifying $\sqrt{\tau}$ as the step size and $\gamma = 2\friction$ as the friction. 
Note the different time scales in~\cref{eq:overdamped_limit_ila,eq:ILA_is_underdamped} which are a consequence of the fact that the underdamped system corresponds to a second order differential equation.
The form~\cref{eq:ILA_is_underdamped} also exposes the structure of the scheme.
Namely, that~\eqref{eq:ila} is a Gauss-Seidel style discretization of underdamped Langevin where the update of each variable uses already the updated version of the other one.\footnote{See also~\Cref{SM:ssec:gibbs} for a relation between~\eqref{eq:ila} and overrelaxed Gibbs sampling.}

To discretize \cref{eq:uld_intro}, there is a vast literature on so-called \emph{splitting} methods
~\cite{leimkuhler2013rational,leimkuhler2024kinetic}.
These separate the right-hand side of~\cref{eq:uld_intro} into individual components that can be integrated exactly,
\begin{equation}
    \dd \begin{bmatrix}
        X_t\\
        V_t
    \end{bmatrix}
    = 
    \underbrace{
    \begin{bmatrix}
        V_t\\
        0
    \end{bmatrix}
    }_{\mathcal{A}}
    \dd t
    +
    \underbrace{
    \begin{bmatrix}
        0\\
        -\nabla U(X_t)
    \end{bmatrix}
    }_{\mathcal{B}}
    \dd t+
    \underbrace{
    \begin{bmatrix}
        0\\
        -\gamma V_t + \sqrt{2\gamma}
    \end{bmatrix}
    }_{\mathcal{O}}
    \dd W_t.
\end{equation}
Composing the Markov kernels induced by the propagators $\mathcal{A}$, $\mathcal{B}$, and $\mathcal{O}$ in different orders yields a variety of schemes (\cf~\cite{leimkuhler2013rational,leimkuhler2024kinetic}) where the maximum number of appearances of any propagator is called the order of the scheme. Second-order splittings with symmetric structure have become particularly popular due to favorable properties such as bias reduction~\cite{leimkuhler2013rational,leimkuhler2024kinetic}. As a representative of these methods in this work we also consider
the scheme obtained from the ordering~\eqref{eq:baoab}, where every propagator appearing twice is ever applied with half the step size. Denoting the step size as $\step$,~\cref{eq:baoab}, thus, reads as
\begin{equation}\label{eq:baoab1}
    \begin{cases}
        V^{k+1/3} =& V^k - \frac{\step}{2}\nabla U(X^k)\\
        X^{k+1/2} =& X^k+\frac{\step}{2}V^{k+1/3}\\
        V^{k+2/3} =& \expo{-\gamma \step}V^{k+1/3} + \sqrt{1-\expo{-2\step\gamma}}Z^k\\
        X^{k+1} =& X^{k+1/2}+\frac{\step}{2}V^{k+2/3}\\
        V^{k+1} =& V^{k+2/3} - \frac{\step}{2}\nabla U(X^{k+1})
    \end{cases}
\end{equation}
Denoting $\beta = \expo{-\gamma\step}$, $\tau = \frac{\step^2}{2}(1+\beta)$ and rewriting~\cref{eq:baoab1} into a single variable scheme yields
\begin{equation}\label{eq:baoab}
    X^{k+1} = X^k - \tau \nabla U (X^k) + \beta(X^k - X^{k-1}) + \sqrt{2\tau (1 - \beta)} \frac{N^{k-1} + N^k}{2}.
    \tag{\ensuremath{\mathcal{BAOAB}}}
\end{equation}
Remarkably, this scheme is~\cref{eq:ila} up to the fact that the noise contains also a contribution of the previous iteration rendering~\cref{eq:baoab} non-Markovian in this form. There is, however, an even more astonishing connection to~\eqref{eq:ila}. Indeed, one can easily verify that if $U$ is quadratic, that is, $\nabla U$ is linear, the~\eqref{eq:baoab} sequence is equivalent to the sequence generated by the averages of two consecutive~\cref{eq:ila} iterates, that is $Z^k = X^k_{\mathrm{BAOAB}}$ where we define 
\begin{equation}\label{eq:debiasing}
    Z^k\coloneqq\frac{X^{k}_\mathrm{ILA} + X^{k-1}_\mathrm{ILA}}{2}.\tag{ILA+}
\end{equation}
We will see below that this averaging leads to a~\emph{de-biasing}, that is, $(Z^k)_k$ admits precisely $\pi$ as its stationary distribution, implying also unbiasedness of~\cref{eq:baoab}.
Inspired by splitting methods, a promising additional extension of~\cref{eq:ila} would be a symmetrized version thereof, where we perform two half steps for the velocity component and inbetween a full step for the position component in~\cref{eq:ILA_is_underdamped}.

By analogy with the optimization setting, it is conjectured that suitable discretizations of~\eqref{eq:uld_intro} should yield samplers with $\tilde\Oc(\sqrt{\kappa})$ complexity for strongly convex targets.
In the continuous-time setting this has indeed been established~\cite{li2026space,lu2020explicit}.
For quadratic potentials, in this work we achieve the desired acceleration, \ie,~\eqref{eq:ila} attains $\tilde\Oc(\sqrt{\kappa})\footnote{The tilde hides logarithmic factors.}$ complexity.
Besides our work, the scaling $\tilde\Oc(\sqrt{\kappa})$ of other kinetic Langevin discretizations has been shown for Gaussians in~\cite{gouraud_hmc_2025,jiang2023dissipation}.
Interestingly, contrary to these works, we achieve $\tilde\Oc(\sqrt{\kappa})$ without resorting to a higher-order splitting.
Moreover, the complexity in~\cite{gouraud_hmc_2025} is defined in an asymptotic manner and scales as $\tilde\Oc(\sqrt{d\kappa})$ in the dimension whereas in our case the dimension enters only logarithmically. Note that the complexity $\tilde\Oc(\sqrt{\kappa})$ is optimal for Gaussians~\cite{chewi2024query}.

Beyond the Gaussian setting, we will see that~\eqref{eq:ila} is a discretization of~\eqref{eq:uld_intro} that remains geometrically ergodic for a wider parameter range than previously previously available in the literature, allowing, in particular, for arbitrarily small friction $\gamma$ and larger step-sizes.

During the last days of writing this manuscript, independent of our work the excellent preprint~\cite{bou2026provable} has been published which uses the analogy between heavy-ball and underdamped Langevin as well.\footnote{the first version of our manuscript had already been online for several months at this point} There it is shown, that the complexity $\tilde\Oc(\sqrt{\kappa})$ in total variation cannot be achieved by~\cref{eq:baoab} and similar splittings uniformly over the class of all strongly convex potentials with condition number smaller than $\kappa$. Thus, to achieve ballistic acceleration, randomized schemes may be needed. The present article is still of interest, as both the proven optimality on Gaussians as well as our analysis in the general (even non-logconcave) setting are unaffected by the negative result in~\cite{bou2026provable}. Moreover, we want to mention that~\cite{bou2026provable} considers only convergence of the schemes to their (biased) invariant distributions and, contrary to us, does not investigate the bias.

We summarize our contributions.
\paragraph{Contributions}
\begin{enumerate}
    \item We introduce~\eqref{eq:ila}, a gradient-based \gls{mcmc} method for accelerated sampling from Gibbs distributions.
    The scheme is simple and poses the direct sampling pendant of the heavy ball method from optimization~\cite{polyak1964some}.
    \item In the Gaussian setting, we show that~\cref{eq:ila} is geometrically ergodic and derive explicit representations of its stationary distribution allowing us to estimate its bias to $\pi$.
    In particular, despite the scheme being only first order, we obtain a bias of order two and with the de-biasing strategy~\cref{eq:debiasing} the scheme becomes unbiased. As a by-product we also show that~\cref{eq:baoab} is unbiased for Gaussian targets.
    As a result we achieve $\delta$-accurate samples in Wasserstein-2 distance with $\tilde\Oc(\sqrt{\kappa})$ iterations.
    \item In the non-Gaussian and even non-convex setting, we show that the scheme is geometrically ergodic under the simple conditions $\beta\in (0,1)$, $\tau<\frac{2(1-\beta)}{L}$, where $L$ is the Lipschitz constant of $\nabla U$.
    In particular, contrary to the literature, this condition allows for arbitrarily small friction coefficients ($\beta \to 1$), which is conjectured to be crucial for acceleration~\cite{altschuler2026shifted}.
    In addition, using \cref{eq:ila}'s relation to underdamped Langevin dynamics, we establish vanishing bias to $\pi$ as $\tau\rightarrow 0$, also in the general case.
    \item Finally, in \Cref{sec:numerical}, we provide numerical experiments confirming the efficacy and efficiency of~\eqref{eq:ila}. The experiments cover ill-conditioned Gaussians, confirming the $\sqrt{\kappa}$ complexity, total variation image denoising and maximum likelihood learning of an energy-based model for molecular structure generation.
\end{enumerate}

\section{Notation, Preliminaries, and Assumptions}
We denote the Borel $\sigma$-algebra on $\R^d$ as $\Bc(\R^d)$ and the space of all probability measures on $\Bc(\R^d)$ as $\Pc(\R^d)$.
The subspace of all probability measures with finite $p$-th moment is denoted as $\Pc_p(\R^d)$. For any $\mu\in \Pc(\R^d)$, and a measurable function $f$ we denote $\mu(f) = \int f\dd \mu$. The Lebesgue measure on $\R^d$ is denoted as $\lambda^d$.
The Wasserstein-$2$ distance is defined as $\Wc_2^2(\mu,\nu) \coloneqq \min_{\gamma\in\Pi(\mu,\nu)}\int_{\R^{2d}} \|x-x'\|^2\dd\gamma(x,x')$ for $\mu,\nu\in\Pc_2(\R^d)$,
where $\Pi(\mu,\nu)$ is the set of all couplings of $\mu$ and $\nu$, that is, the set of all probability measures $\gamma\in\Pc(\R^{2d})$ with marginals $\mu$ and $\nu$~\cite{villani2008optimal}.
We refer to $\hat{\gamma}\in\Pi(\mu,\nu)$ realizing the minimum in $\Wc_2(\mu,\nu)$ as well as a tuple of random variables $(X,Y)\sim \hat{\gamma}$ as an \emph{optimal coupling} of $(\mu,\nu)$. We will use the notation $\Wc_2(X,Y) \coloneqq \Wc_2(\law(X),\law(Y))$ for random variables $X,Y$.
For measurable $\Vc:\R^{2d}\rightarrow [0,\infty]$ we define $\Vc$-norm $\|\mu\|_\Vc = \sup_{f\leq \Vc} \int f\dd \mu$.
A Markov kernel is a mapping $M: \R^d\times \Bc(\R^d)\rightarrow [0,1]$ such that $M(x,\cdot)$ is a probability measure for any $x\in\R^d$ and $M(\cdot,A)$ is measurable for any $A\in\Bc(\R^d)$.
We denote the application of a Markov kernel $M$ to a probability measure $\mu\in\Pc(\R^d)$ as $\mu M(A) = \int_{\R^d} M(x,A)\dd \mu(x),\quad A\in\Bc(\R^d)$.
A Markov kernel can also be viewed as a linear operator on bounded measurable functions via the adjoint operation of of $\mu\mapsto \mu M$.
A Markov semigroup is a family \((M_t)_{t \geq 0}\) of such operators satisfying \(M_t \1 = \1\) with $\1$ the constant function with value one, \(M_t f \geq 0\) whenever \(f\geq 0\), \(M_t f \to f\) in \(L^2(\mathbb{R}^d,\mu)\) as \(t \to 0\) where $\mu$ is the invariant measure for $M_t$, and the semigroup property \(M_t \circ M_s = M_{t+s}\)~\cite{bakry2013analysis}.
For the kinetic Langevin dynamics~\cref{eq:uld_intro}, the associated semigroup is denoted \((P_t)_{t \geq 0}\). 
That is, if \((X_t, V_t)\) solves the dynamics with initial law \(\mu = \law(X_0, V_0)\), then \(\mu P_t = \law(X_t, V_t)\), and we often write \(Z_t = (X_t, V_t)\) to denote position-velocity pairs.

\section{Main Results}
In this section we prove our main theoretical results, first for Gaussian targets and then for general potentials.
\subsection{The Gaussian Case}
For Gaussian targets we recover analogues of the classical convergence results for the discrete heavy-ball method, with one crucial difference: unlike in optimization, in sampling we may also be interested in the regime $\tau \to 0$, if the scheme is biased (\cf, however,~\Cref{prop:debiases}).
Throughout, $\kappa = L/m$  denotes the condition number, with $L$ and $m$ the largest and smallest eigenvalues of $\nabla^2 U$, respectively.
\begin{theorem}
    [Analysis in the Gaussian case]\label{thm_gauss}
    Assume $U(x) = \frac{1}{2}(x-p)^\top\Lambda (x-p)$ with $\Lambda\in\R^{d\times d}$ symmetric, with $0\prec m\preceq \Lambda\preceq L$, and $p\in \R^d$.
    We have the following:
    \begin{enumerate}
        \item \emph{Ergodicity:} For any $\beta\in [0,1)$ and $\tau\in (0,2(1+\beta)/L)$, and any initial distribution $\mu \in \mathcal P_2(\R^d)$,~\eqref{eq:ila} is geometrically ergodic in the Wasserstein-2 distance with its unique stationary distribution denoted as $\pi_{\tau,\beta}$. That is, there exists $\rho=\rho(\tau,\beta)\in (0,1)$ such that for any $\delta\in (0,1-\rho)$ there exists $K = \tilde\Oc(\delta^{-1})$ such that for $k\geq K$
        \begin{equation}\label{eq:gauss_proof_convergence}
            \Wc_2(X^k,\pi_{\tau,\beta})\leq (\rho+\delta)^{k}\bigg(\E\left[ \|X^0-p\|^2 + \|X^1-p\|^2 \right]+ \frac{2\tau(1-\beta)d}{1-(\rho+\delta)^2}\bigg)^{1/2}.
        \end{equation}
        \item \emph{Optimal $\beta$:} For every admissible $\tau$ the optimal convergence rate, \ie, the minimal $\rho$ is achieved for
        \begin{equation}\label{eq:beta-optimal}
            \beta^\ast(\tau) \coloneqq \argmin_{\beta\in \left(\frac{L\tau}{2}-1,1\right)} \rho(\tau,\beta) = (1-\sqrt{L\tau})^2\vee (1-\sqrt{m\tau})^2
        \end{equation}
        with rate $\rho(\tau,\beta^*(\tau))=\sqrt{\beta^*(\tau)}$. 
        The minimal $\rho$ over all admissible $\tau$ \emph{and} $\beta$ is
        \begin{equation}\label{eq:optimal_tau_beta}
            \rho = \frac{1-1/\sqrt{\kappa}}{1+1/\sqrt{\kappa}}\quad\text{for}\quad \tau = \frac{4}{(\sqrt{L}+\sqrt{m})^2} \quad \text{and} \quad \beta = \left(\frac{1-1/\sqrt{\kappa}}{1+1/\sqrt{\kappa}}\right)^2.
        \end{equation}
        \item \emph{Bias:} If $\tau\leq \frac{1}{L}$ it holds $\Wc_2(\pi_{\tau,\beta},\pi)\leq \frac{\tau}{2}\sqrt{\trace(\Lambda)}$.
    \end{enumerate}
\end{theorem}
\begin{remark}
    Some remarks are in order.
    \begin{enumerate}
        \item Bias of order $\Oc(\tau)$ is achieved also under the condition $\tau<\frac{2(1+\beta)}{L}$ but the additional assumption $\tau\leq 1/L$ allows for a simpler bound.
        \item Above we recover the same optimal convergence rate as for the classical heavy-ball method for optimization. However, since contrary to optimization, in sampling there is also an incentive to decrease the step size in order to reduce the bias of the stationary distribution, we also derive for every step size $\tau$ the optimal $\beta$. Below we show that~\cref{eq:debiasing} becomes unbiased for any admissible step size.
        \item As mentioned above (\cf~\eqref{eq:ILA_is_underdamped}), under the optimal value $\beta^*(\tau)$,~\eqref{eq:ila} poses a discretization of the kinetic Langevin dynamics with step size $\sqrt{\tau}$. 
        \item This results unifies the convergence in both regimes, fixed $\beta$, where~\cref{eq:ila} corresponds to~\cref{eq:OD} and $1-\beta = \Oc(\sqrt{\tau})$ where~\cref{eq:ila} corresponds to~\cref{eq:uld_intro}. Considering the different time scales $\tau$ vs $\sqrt{\tau}$ we obtain order 1 and order 2 bias, respectively. The latter is achieved without relying on a higher-order splitting.
        \item Note that our optimal choice of $\beta^*$ directly leads to an optimal friction parameter choice in the underdamped Langevin dynamics~\eqref{eq:uld_intro}, namely, $\gamma = 2\sqrt{m}$.
    \end{enumerate}
\end{remark}
\begin{proof}
    Note that, without loss of generality, we can assume $p=0$ and a diagonal precision matrix $\Lambda = \diag(\lambda_1,\dots,\lambda_d)$. 
    Rewriting the recursion as a system and denoting $\sigma = \sqrt{2\tau(1-\beta)}$ yields
    \begin{equation}\label{eq:ila-recursion}
        Y^{k+1}\coloneqq
        \begin{bmatrix}
            X^{k+1}\\
            X^k
        \end{bmatrix}
        =
        \underbrace{\begin{bmatrix}
            1-\Lambda\tau + \beta&-\beta\\
            1&0
        \end{bmatrix}}_{\eqqcolon Q}
        \begin{bmatrix}
            X^{k}\\
            X^{k-1}
        \end{bmatrix}
        +
        \begin{bmatrix}
            \sigma N^k\\
            0
        \end{bmatrix}.
    \end{equation}
    We may unroll the iteration as
    \begin{equation}\label{eq:heavy_ball_unrolled}
        Y^k 
        = Q^k Y^0 + \sum_{\ell=0}^{k-1} Q^{\ell}
        \begin{bmatrix}
            \sigma N^{k-1-\ell}\\
            0
        \end{bmatrix}
        \overset{d}{=}
        Q^k Y^0 + \sum_{\ell=0}^{k-1} Q^{\ell}
        \begin{bmatrix}
            \sigma N^{\ell}\\
            0
        \end{bmatrix}.
    \end{equation}
    In what follows, we work with the version on the right-hand side, where we re-indexed the Gaussian \iid~sequence $(N^k)_k$.
    $Y^k$ is a Gaussian with mean $\E[Y^k] = Q^k\E[Y^0]$ and, due to independence, the covariance reads as
    \begin{equation}\label{eq:cov_y}
        \cov[Y^k] = Q^k\cov[Y^0](Q^\top)^k + \sigma^2\sum_{\ell=0}^{k-1}Q^\ell\begin{bmatrix}
            \id&0\\
            0&0
        \end{bmatrix}(Q^\top)^\ell.
    \end{equation}
    Ergodicity, thus, follows if the spectral radius of $Q$, $\rho(Q) \coloneqq \max_i |\eta_i|$ where $\eta_i$ are the eigenvalues of $Q$,
    is strictly smaller than one. The proof is inspired by the classical one from optimization. However, we repeat it here for completeness and because dealing with the variance of $Y^k$ is not covered by the classical result.
    
    Let $(u,v)\in\R^{2d}$ be an eigenvector of $Q$ with corresponding eigenvalue $\eta$. Due to the block diagonal structure of $Q$, this implies that for some $\lambda\in \{\lambda_i\;|\; i=1,\dots,d\}$ we have $\eta u = (1-\lambda\tau + \beta)u - \beta v$, $\eta v = u$.
    It follows immediately that neither $u$ nor $v$ can be zero, as otherwise $u=v=0$.
    Inserting the second into the first equation leads to $\eta^2-(1-\lambda\tau + \beta)\eta + \beta  = 0$.
    One can check that for $0\leq \beta <1$ and $0<\tau <\frac{2(1+\beta)}{L}$, it follows that $\rho(Q)<1$ (\cf~\cite[Theorem 1, Section 3.2.1]{polyak1987introduction}).
    As a consequence, for any $\delta\in (0,1-\rho(Q))$, there exists $K=\tilde\Oc(\delta^{-1})$ such that for $k\geq K$, $\|Q^k\|\leq (\rho(Q) + \delta)^k$~(\cf~\Cref{SM:convergence_corollary} for details regarding the scaling of $K$). In particular, we have by independence for any $K\leq k_1<k_2$
    \begin{equation}
        \begin{aligned}
            \E \left[ \left\|\sum_{\ell=k_1}^{k_2}Q^\ell \begin{bmatrix}
            \sigma N^\ell\\0
            \end{bmatrix} \right\|^2\right]
            =\E \left[ \sum_{\ell=k_1}^{k_2} \left\|Q^\ell \begin{bmatrix}
                \sigma N^\ell\\0
            \end{bmatrix} \right\|^2\right]
            \leq& \sigma^2 d \sum_{\ell=k_1}^{k_2} (\rho(Q)+\delta)^{2\ell},
        \end{aligned}
    \end{equation}
    implying that $(Y^k)_k$, or more accurately $(\law(Y^k))_k$, is a Cauchy sequence in $(\Pc_2(\R^{2d}),\Wc_2)$\footnote{Note that the re-indexing of the $(N^\ell)_\ell$ prohibits the conclusion of a Cauchy sequence in $L^2(\P)$.}.
    By completeness, $(\law(Y^k))_k$, thus, converges to its unique stationary distribution whose first marginal\footnote{Indeed, at stationarity, both marginals are identical.} we denote as $\pi_{\tau,\beta}$. It follows for $k\geq K$
    \begin{equation}\label{eq:gaussian}
        \begin{aligned}
            \Wc_2^2(X^k,\pi_{\tau,\beta})
            \leq& \E\left[ \left\|Q^k Y^0\right\|^2 \right]
                +\E\left[ \left\|\sum_{\ell=k}^{\infty}Q^\ell \begin{bmatrix}
                \sigma N^\ell\\0
                \end{bmatrix} \right\|^2 \right]\\
            \leq& (\rho(Q)+\delta)^{2k}\E\left[ \|Y^0\|^2 \right]+\sigma^2d  \sum_{\ell=k}^{\infty}(\rho(Q)+\delta)^{2\ell}\\
            \leq& (\rho(Q)+\delta)^{2k}\bigg(\E\left[ \|X^0\|^2 + \|X^1\|^2 \right]+ \frac{\sigma^2d}{1-(\rho(Q)+\delta)^2}\bigg),
        \end{aligned}
    \end{equation}
    concluding the proof of~\eqref{eq:gauss_proof_convergence} with $\rho(\tau, \beta) \coloneqq \rho(Q(\tau, \beta))$ in a slight abuse of notation.
    For the optimal $(\tau,\beta)$, see the supplement,~\Cref{SM:lemma:supp:optimal_beta}.
    
    Next, let us investigate the stationary distribution and its bias to the target. Note, moreover, that due to $\Lambda$ being diagonal all four blocks of $Q^\ell$ are diagonal as well and we may consider without loss of generality the one-dimensional case, \ie, $\Lambda = \lambda_i\in\R$ where we drop the index $i$ for simplicity and just write $\lambda$.
    By L\'evy's continuity, $\pi_{\tau,\beta}$ is Gaussian again so that its mean and covariance fully characterize $\pi_{\tau,\beta}$.
    From above, we have that its mean is zero.
    Recall that $\cov[X^{k}]$ is the upper left entry of $\cov[Y^k]$. One can easily check that
    \begin{equation}
        \begin{bmatrix}
        Q^\ell\begin{bmatrix}
            \id&0\\
            0&0
        \end{bmatrix}
    (Q^\top)^\ell\end{bmatrix}_{1,1}
    = (q^{\ell}_{1,1})^2
    \end{equation}
    where $q^\ell_{i,j}$ denotes the entry $(i,j)$ of $Q^\ell$, not to be confused with $(q_{i,j})^\ell$.
    We find the recursion $q^{k+1}_{1,1} = q_{1,1} q^{k}_{1,1}-\beta q^{k-1}_{1,1}$ with $q_{1,1} = 1-\lambda\tau+\beta$. Using the ansatz $q^k_{1,1} = c \eta^k$ we find $q^k_{1,1} = c_1 \eta_+^k + c_2\eta_-^k$
    where $\eta_\pm$ are the solutions of the characteristic equation $\eta^2-a_1\eta +\beta = 0$ and, thus, identical to the eigenvalues of $Q$, and $c_1, c_2$ are determined by the initial conditions $q^0_{1,1} = 1$, $q^1_{1,1} = 1-\lambda\tau + \beta$. 
    \paragraph{Case 1: Fixed $\beta$}
    Assume for simplicity the characteristic equation admits two distinct solutions $\eta_+\neq \eta_-$. The converse case will be handled below. Note, however, that $\eta_+= \eta_-$ is true if and only if $\tau = (1\pm\sqrt{\beta})^2/\lambda$ which can be avoided by simply decreasing $\tau$ by an arbitrarily small amount in the case of $\beta$ being fixed.
    Then we have
    \begin{equation}
        q^k_{1,1} = \frac{\eta_+^{k+1} - \eta^{k+1}_-}{\eta_+-\eta_-}.
    \end{equation}
    As $|\eta_\pm|\leq \rho(Q)<1$, by the formula for the geometric series we have
    \begin{equation}\label{eq:gaussian_var}
        \begin{aligned}
            \cov[\pi_{\tau,\beta}]_{ii}=\lim_{k\rightarrow\infty }\cov[X^k]_{ii} &= \sigma^2\sum_{\ell\geq 0}(q^\ell_{1,1})^2\\
            &= \frac{\sigma^2}{(\eta_+-\eta_-)^2}\left\{ \frac{\eta_+^2}{1-\eta_+^2} + \frac{\eta_-^2}{1-\eta_-^2} - 2\frac{\eta_+\eta_-}{1-\eta_+\eta_-}\right\}\\
            &= \frac{\sigma^2}{(\eta_+-\eta_-)^2}\left\{ \frac{1}{1-\eta_+^2} + \frac{1}{1-\eta_-^2} - \frac{2}{1-\eta_+\eta_-}\right\}\\
            &=  \frac{\sigma^2(1+\eta_+\eta_-)}{(1-\eta_+^2)(1-\eta_-^2)(1-\eta_+\eta_-)}.
        \end{aligned}
    \end{equation}
    Note that $\eta_+\eta_- = \beta$ so that we obtain
    \begin{equation}\label{eq:variance}
        \cov[\pi_{\tau,\beta}]_{ii}=
        \frac{\sigma^2(1+\beta)}{(1-\eta_+^2)(1-\eta_-^2)(1-\beta)}.
    \end{equation}
    By definition of $\eta_\pm$ as roots of $\eta^2 - q_{1,1}\eta + \beta$ we may write $\eta^2 - q_{1,1}\eta + \beta = (\eta_+-\eta)(\eta_--\eta)$, so that by once inserting $\eta=1$ and once $\eta = -1$ we find $(1-\eta_+^2)(1-\eta_-^2)=\lambda\tau(2-\lambda\tau+2\beta)$.
    Inserting into~\eqref{eq:variance} and noting that $\sigma^2 = 2(1-\beta)\tau$ we find
    \begin{equation}
        \begin{aligned}
            \cov[\pi_{\tau,\beta}]_{ii}
            =\frac{2\tau(1-\beta)(1+\beta)}{\lambda\tau(2-\lambda\tau+2\beta)(1-\beta)}
            &=\frac{2(1+\beta)}{\lambda(2-\lambda\tau+2\beta)}\\
            &=\frac{1}{\lambda}\frac{1}{1-\frac{\lambda\tau}{2(1+\beta)}}\\
            &=\frac{1}{\lambda}\sum_{k=0}^\infty\left(\frac{\lambda\tau}{2(1+\beta)}\right)^k\\
            &=\frac{1}{\lambda} + \frac{\tau}{2(1+\beta)}\sum_{k=1}^\infty\left(\frac{\lambda\tau}{2(1+\beta)}\right)^{k-1}\\
            &\leq\frac{1}{\lambda} + \frac{\tau}{1+2\beta},
        \end{aligned}
    \end{equation}
    where the last inequality follows by the formula for the geometric series under the assumption $\tau\leq \frac{1}{L}$.
    Noting that $\cov[\pi_{\tau,\beta}]_{ii}\geq 1/\lambda=\cov[\pi]_{ii}$, and using the formula for the Wasserstein distance between Gaussian measures together with $\sqrt{a+\tau b}\leq\sqrt{a}+\tau b/(2\sqrt{a})$ for all $a,b>0$ and $\tau>0$ by concavity, we obtain, upon summing over all eigenvalues $\lambda_i$,
    \begin{equation}
        \begin{aligned}
            \Wc_2^2(\pi_{\tau,\beta},\pi) 
            \leq \sum_{i}\bigg(\sqrt{\frac{1}{\lambda_i} + \tau} - \frac{1}{\sqrt{\lambda_i}}\bigg)^2
            \leq \sum_{i} \frac{\tau^2 \lambda_i}{4}
            =& \frac{\tau^2}{4}\trace(\Lambda).
        \end{aligned}
    \end{equation}
    \paragraph{Case 2: Variable $\beta$}
    Next, we assume $\beta = \beta(\tau) = (1-\friction\sqrt{\tau})^2$, which is inspired by the form of the optimal friction derived above. The case $\eta_+\neq\eta_-$ is treated exactly as in the fixed $\beta$ setting. 
    Now, however, the converse case $\eta_+ = \eta_-$ is more relevant, as it will strike whenever $\epsilon^2 = \lambda$ and, thus, cannot necessarily be avoided by decreasing the step size.
    When $\eta_+ = \eta_-$, we have $q^k_{1,1} = (1+k)\beta^{k/2}$.
    It follows that
    \begin{equation}\label{eq:degenerat_eta}
        \begin{aligned}
            \cov[\pi_{\tau,\beta}]_{ii}
            =\sigma^2\sum_{\ell\geq 0}(q^\ell_{1,1})^2
            = \frac{\sigma^2(1+\beta)}{(1-\beta)^3}
            = \frac{2\tau(1+\beta)}{(1-\beta)^2}
            &= \frac{2(2-2\friction\sqrt{\tau} + \friction^2\tau)}{\friction^2(2 - \friction\sqrt{\tau})^2}\\
            &= \frac{1}{\lambda} + \frac{\tau}{(2-\sqrt{\lambda}\sqrt{\tau})^2}\\
            &\leq \frac{1}{\lambda} + \tau,
        \end{aligned}
    \end{equation}
    where the last inequality follows once again from $\tau\leq 1/{L}$.
    In total we, therefore, still find
    \begin{equation}
        \begin{aligned}
            \Wc_2^2(\pi_{\tau,\beta},\pi) 
            \leq \sum_{i} \bigg(\sqrt{\frac{1}{\lambda_i} + \tau} - \frac{1}{\sqrt{\lambda_i}}\bigg)^2
            \leq \sum_{i} \frac{\tau^2\lambda_i}{4}
            =& \frac{\tau^2}{4}\trace(\Lambda).
        \end{aligned}
    \end{equation}
\end{proof}
\begin{corollary}[Complexity]\label{cor:ila-complexity}
    Assume~\eqref{eq:ila} is initialized at the mean of the Gaussian target. Set $\beta = (1-\sqrt{m\tau})^2$. In order to achieve $\Wc_2(X^k,\pi)<\delta$,~\eqref{eq:ila} requires $\tilde\Oc(\frac{\sqrt{\kappa}\vee(\frac{d^{1/4}L^{1/4}}{\sqrt{m}})}{\sqrt\delta})$ steps where $\tilde\Oc$ hides logarithmic factors.
\end{corollary}
\begin{remark}\
    \begin{itemize}
        \item Initializing at the mean is not necessary but renders the terms simpler.
        Any initialization with $\E[\|X^0-p\|^2 + \|X^1-p\|^2] = \Oc(\poly(d))$ with $p$ the mean leads to the same complexity.
        \item Note that the result shows two regimes: For $d>L$ the dependence on $\kappa$ improves at the cost of dependence on the dimension.
    \end{itemize}
\end{remark}
As mentioned, we can even debias~\cref{eq:ila} by computing averages of ever two consecutive iterates according to~\cref{eq:debiasing}.
\begin{proposition}[Unbiased~\cref{eq:ila}]\label{prop:debiases}
    Let $Z^k$ be generated according to~\cref{eq:debiasing}, \ie, $Z^k = (X^k+X^{k-1})/2$ with $(X^k)_k$ following~\cref{eq:ila}.
    Then for any $\beta\in (0,1)$, $\tau<\frac{2(1+\beta)}{L}$ it holds
    \begin{equation}
        \begin{aligned}
            \Wc_2(Z^k,\pi)
            \leq& (q+\delta)^{k-1}\bigg(\E\left[ \|X^0-p\|^2 + \|X^1-p\|^2 \right]+ \frac{2\tau(1-\beta)d}{1-(q+\delta)^2}\bigg)^{1/2}
        \end{aligned}
    \end{equation}
    That is, $(Z^k)^k$ converges to the target $\pi$ without any bias. As a consequence, since~\cref{eq:debiasing} is equivalent to~\cref{eq:baoab}, also the latter is unbiased and converges at the same rate.
\end{proposition}
\begin{proof}
    Note that in~\Cref{thm_gauss} we have, in fact, shown that the sequence $Y^k = (X^k,X^{k-1})^\top$ is Cauchy with respect to $\Wc_2$. Since $Z^k = \frac{1}{2} (1,1)Y^k$ we can conclude that also $(Z^k)_k$ is Cauchy with exactly the same convergence rate.
    Moreover, we have
    \begin{equation}
        \begin{aligned}
            \lim_{k\rightarrow\infty}\cov[Z^k]
            =&  \frac{1}{4}\lim_{k\rightarrow \infty}\cov[X^k] + \cov[X^{k-1}] + 2\cov[X^k,X^{k-1}]\\
            =& \frac{1}{2}\lim_{k\rightarrow \infty}\cov[X^k] + \cov[X^k,X^{k-1}].
        \end{aligned}
    \end{equation}
    We have already derived $\lim_{k}\cov[X^k]$ above so that it remains to derive the \emph{mixed} covariance.
    Note that we can again restrict to the 1D case, as different coordinates of the sequence $(X^k)_k$ are independent.
    Since $\cov[X^k]$ is the off-diagonal entry of~$\cov[Y^k]$ from~\cref{eq:cov_y}, we require
    \begin{equation}
        \begin{bmatrix}
            Q^k \begin{bmatrix}
                \id&0\\
                0&0
            \end{bmatrix}(Q^k)^\top
        \end{bmatrix}_{1,2} = q^k_{1,1}q^k_{2,1}.
    \end{equation}
    Moreover, it holds that $q^k_{2,1} = q^{k-1}_{1,1}$. Due to~\cref{eq:cov_y} and the explicit form of $q^k_{1,1}$ we then find
    \begin{equation}
        \begin{aligned}
            \lim_{k\rightarrow \infty}\cov[X^k,X^{k-1}] 
            =& \sigma^2\sum_{k=0}^\infty q^k_{1,1}q^{k-1}_{1,1}\\
            =& \frac{\sigma^2}{(\eta_+-\eta_-)^2}\sum_{k=0}^\infty (\eta_+^{k+1} - \eta_-^{k+1})(\eta_+^{k} - \eta_-^{k})\\
            =&\frac{\sigma^2}{(\eta_+-\eta_-)^2}\left\{\frac{\eta_+}{1-\eta_+^2} + \frac{\eta_-}{1-\eta_-^2}- \frac{\eta_+ + \eta_-}{1-\eta_+\eta_-}\right\}\\
            =&\sigma^2\frac{\eta_++\eta_-}{(1-\eta_+^2)(1-\eta_-^2)(1-\eta_+\eta_-)}.
        \end{aligned}
    \end{equation}
    Using as in the proof of~\Cref{thm_gauss} that $\eta_+\eta_- = \beta$, $(1-\eta_+^2)(1-\eta_-^2) = \lambda\tau(2+2\beta-\lambda\tau)$ and $\eta_++\eta_- = 1-\lambda\tau + \beta$ we find
    \begin{equation}
        \begin{aligned}
            \lim_{k\rightarrow \infty}\cov[X^k,X^{k-1}] 
            =& \sigma^2\frac{1-\lambda\tau + \beta}{\lambda\tau(2+2\beta-\lambda\tau)(1-\beta)}.
        \end{aligned}
    \end{equation}
    Combining with the results from~\Cref{thm_gauss} we find
    \begin{equation}
        \begin{aligned}
            \lim_{k\rightarrow\infty}\cov[Z^k]
            = \frac{\sigma^2}{2}\frac{(1+\beta) + (1-\lambda\tau + \beta)}{\lambda\tau(2+2\beta-\lambda\tau)(1-\beta)}
            =& \frac{\sigma^2}{2\lambda\tau(1-\beta)} = \frac{1}{\lambda}.
        \end{aligned}
    \end{equation}
\end{proof}
Since~\cref{eq:debiasing} is unbiased, we do not have to choose $\tau$ small to achieve $\delta$-accuracy. Instead, we can simply choose the optimal momentum and step size~\cref{eq:optimal_tau_beta} and obtain the following complexity where we get rid of polynomial dependence on $\delta$.
\begin{corollary}[Complexity of~\cref{eq:debiasing}]\label{cor:complexity_ila+}
    \eqref{eq:debiasing}, and therefore, also~\eqref{eq:baoab} with the optimal momentum and step size~\cref{eq:optimal_tau_beta} require $\tilde\Oc(\sqrt{\kappa})$ iterations to achieve $\Wc_2(Z^k,\pi)<\delta$.
\end{corollary}

\begin{figure}
    \centering
    \includegraphics[]{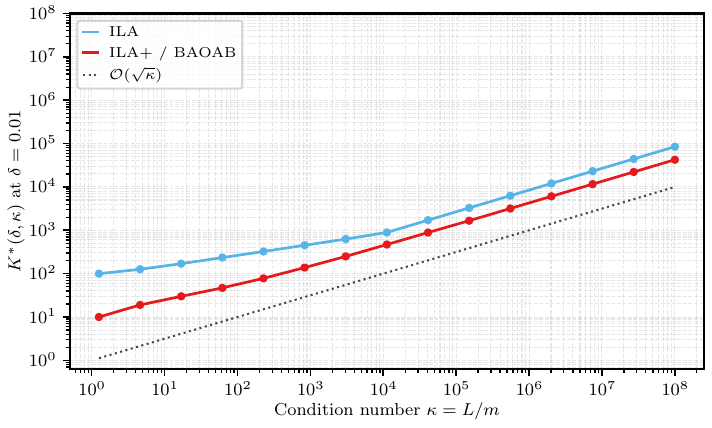}
    \caption{Empirical verification of the complexity of~\eqref{eq:ila} and~\cref{eq:debiasing} (\cf~\Cref{cor:ila-complexity,cor:complexity_ila+})
    on a bivariate Gaussian potential; 
    For~\cref{eq:debiasing} we choose $\tau$, $\beta$ as the optimal parameters~\cref{eq:optimal_tau_beta} and obtain improved complexity.
    }
    \label{fig:ila-complexity}
\end{figure}

\begin{figure}
    \centering
    \includegraphics[]{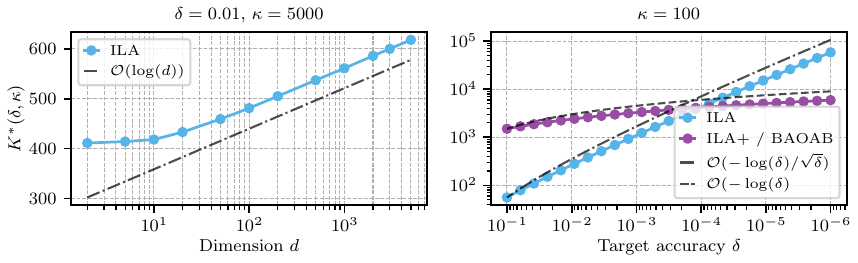}
    \caption{Visualization of $K^\ast (\delta, \kappa) = \min \{k \mid \Wc_2(\pi_k, \pi) \leq \delta\}$ as a function of the dimension $d$ (left) and the target accuracy $\delta$ (right) in comparison to the predictied results for a multivariate Gaussian.}
    \label{fig:complexity-dimension-accuracy}
\end{figure}

\Cref{fig:ila-complexity} confirms the $\tilde\Oc(\sqrt{\kappa})$ scaling of \eqref{eq:ila} and~\cref{eq:baoab} numerically.
In \Cref{fig:complexity-dimension-accuracy}, we also ablate over the $d$ of the sampling problem and the target accuracy $\delta$ confirming the proven complexity. Note the improved scaling with $\delta$ due to the debiasing in~\cref{eq:debiasing}.


\subsection{The General Case}\label{sec:non_gauss}
We consider now non-Gaussian and even non-log-concave targets and impose instead the following conditions on $U$.
\begin{assumption}[General assumptions on $U$]\label[assumption]{ass:assumptions}\
    \begin{enumerate}
        \item $U:\R^d\rightarrow \R$ is twice continuously differentiable.
        \item $U$ is $L$-smooth 
        ,\ie, $\|\nabla^2 U(x)\|\leq L$ for all $x \in \R^d$ and some $L>0$. \label{ass:lipschitz}
        \item There exists $\grc\geq 1$ such that $\liminf_{\|x\|\rightarrow \infty}\inner{x}{\nabla U(x)}\|x\|^{-\grc{}}>0$.\label{ass:growth}
        \item $U(0) = \min_{x} U(x)=0$.
    \end{enumerate}
\end{assumption}
The last assumption is just a convention and can be enforced by shifting $U$.

Motivated by the optimal momentum $\beta^\ast(\tau)$ identified in the Gaussian analysis, we will at times additionally impose the following.
\begin{assumption}[Kinetic momentum]\label[assumption]{ass:beta}
    The parameter $\beta$ is of the form $\beta(\tau) = 1-2\friction\sqrt{\tau}$ for some $\friction>0$.
\end{assumption}
All results below carry over to any $\beta (\tau) = 1 - 2 \varepsilon \sqrt{\tau} + \mathcal O (\tau)$ -- in particular to the choice $\beta(\tau) = (1 - \varepsilon \sqrt{\tau})^2$ from the Gaussian setting -- but we restrict ourself to the above choice for simplicity.

\subsubsection{Discrete Time}
To analyze \cref{eq:ila} theoretically, we first rewrite it in terms of a first-order Markov chain.
With the auxiliary variable $V^{k} = (X^k-X^{k-1})/\sqrt \tau$, \eqref{eq:ila} becomes
\begin{equation}\label{eq:ila-2}
    \begin{cases}
        V^{k+1} &= \beta V^k - \sqrt{\tau}\nabla U(X^k) + \sqrt{2(1 - \beta)}N^k,\\
        X^{k+1} &= X^k +\sqrt{\tau} V^{k+1},
    \end{cases}
    \tag{ILA-II}
\end{equation}
and we denote its Markov transition kernel with $R_{\tau,\beta}$.

Our analysis rests on the classical Lyapunov theory for Markov chains~\cite{meyn2012markov}: geometric ergodicity follows once we have established irreducibility, strong aperiodicity and a Lyapunov drift condition, that is, the existence of a function $\Vc: \R^{2d}\rightarrow [1,\infty)$ unbounded off of small sets such that the Markov kernel $R_{\tau,\beta}$ satisfies for some $\lambda\in (0,1)$, $b>0$
\begin{equation}\label{eq_defin_lyap}
    R_{\tau,\beta}\Vc(x,v)\leq \lambda \Vc(x,v) + b \tag{Lyap}
\end{equation}
Such a condition directly implies boundedness of $R^k_{\tau,\beta}\Vc(x,v)$ for all $\Vc$, as iterating the condition yields $R^k_{\tau,\beta}\Vc(x,v)\leq \lambda^k \Vc(x,v) + \frac{b}{1-\lambda}$. 
Note that the second term explodes as $1-\lambda\to 0$ if we do not ensure $b\to 0$ at the same rate. 
While in the purely discrete setting we may ignore this, if we are interested in the regime $\tau\to 0$ the behavior for $\lambda\to 1$ will be relevant.

The consequence of such a Lyapunov condition together with the other mentioned conditions is the existence of $C>0$, $\kappa\in (0,1)$ such that~\cite[Theorem 15.0.2]{meyn2012markov}
\begin{equation}
    \|\mu R^k_{\tau,\beta} - \pi_{\tau,\beta}\|_\Vc\leq C\kappa^k \mu(\Vc),\quad \mu \in \Pc(\R^{2d}).
\end{equation}
We will focus now on deriving such Lyapunov conditions and discuss the remaining conditions to obtain ergodicity afterward.
In the following we will first derive a condition for large step sizes, ensuring ergodicity up to $\tau<\frac{2}{L}$. The obtained condition will however deteriorate as $\tau\to0$, which is why we will afterward derive a second Lyapunov condition which is only valid for $\tau<\frac{1-\beta^2}{L}$ but also remains stable as $\tau \to 0$.
\paragraph{Lyapunov condition for large step sizes}
The main part of the proof is the following lemma.
\begin{lemma}\label{lemma:omega_stronger}
    Let~\Cref{ass:assumptions} hold. Define $\omega(x,v) = U(x) + \frac{\beta}{2}\|v\|^2$ and denote $q \coloneqq \beta v - \sqrt\tau \nabla U(x)$. Then for any $\beta\in (0,1)$ and $\tau<\frac{2(1-\beta)}{L}$
    \begin{equation}
        R_{\tau,\beta}\omega(x,v)\leq \omega(x,v) - \frac{\beta}{2}\|q-v\|^2 - \|q\|^2\left\{ 1 - \beta - \frac{L\tau}{2}\right\} + 2d(1+\frac{L\tau}{2})(1-\beta).
    \end{equation}
\end{lemma}
\begin{proof}
    Note that with the notation $Q^k = \beta V^k - \sqrt\tau \nabla U(X^k)$ we can write~\eqref{eq:ila} as 
    \begin{equation}
        \begin{cases}
            V^{k+1} =& Q^k + \sqrt{2(1-\beta)} N^k\\
            X^{k+1} =& X^k + \sqrt{\tau} Q^{k} + \sqrt{2\tau (1-\beta)}N^k.
        \end{cases}
    \end{equation}
    By $L$-smoothness of $U$ (\cf~\Cref{ass:assumptions}) and the fact that $\sqrt\tau \nabla U(x) = \beta v - q$ we have
    \begin{equation}
        \begin{aligned}
            R_{\tau,\beta}\omega(x,v)
            \leq& U(x) + \langle\nabla U(x), \sqrt{\tau}q\rangle + \frac{L\tau}{2}\|q\|^2 + \frac{2Ld\tau(1-\beta)}{2}+ \frac{\beta}{2}\|q\|^2 + 2d(1-\beta)\\
            =& U(x) + \beta\langle v, q\rangle  - \|q\|^2+ \frac{L\tau}{2}\|q\|^2 + \frac{\beta}{2}\|q\|^2 + 2d(1+\frac{L\tau}{2})(1-\beta)\\
            =& \omega(x,v) + \frac{\beta}{2}\|q\| - \frac{\beta}{2}\|q-v\|^2 - \|q\|^2+ \frac{L\tau}{2}\|q\|^2 \\
            &+ \frac{\beta}{2}\|q\|^2+ 2d(1+\frac{L\tau}{2})(1-\beta)\\
            =& \omega(x,v) - \frac{\beta}{2}\|q-v\|^2 - \|q\|^2\left\{1 - \beta - \frac{L\tau}{2}\right\} + 2d(1+\frac{L\tau}{2})(1-\beta).
        \end{aligned}
    \end{equation}
\end{proof}
\begin{proposition}[Drift condition for large steps]\label{prop:lyapunov_stronger}
    Let~\Cref{ass:assumptions} hold and assume in addition that $\|\nabla U\|^2\gtrsim U$. Define the Lyapunov function $\Vc(x,v) = \expo{\eta\sqrt{1+\omega(x,v)}}$ with $\omega(x,v)=U(x)+\frac{\beta}{2}\|v\|^2$. Then for every $\beta\in (0,1)$, $\tau<\frac{2(1-\beta)}{L}$ there exist $\lambda\in (0,1)$, $b>0$ such that~\cref{eq_defin_lyap} holds.
\end{proposition}
\begin{proof}
    The proof is inspired by~\cite{durmusUniformMinorizationCondition2025}. Fix $(x,v)\in \R^{2d}$ and denote $\bar v(n) = \beta v - \sqrt{\tau} \nabla U(x) + \sqrt{2(1 - \beta)}n$,
    \ie, the update of the velocity component in~\eqref{eq:ila-2} as a function of the noise $n$.
    One can easily check, \eg, by computing the first derivative, that under~\Cref{ass:assumptions} and $\|\nabla U\|^2\gtrsim U$
    \begin{equation}
        n \mapsto \sqrt{1+\omega\left(x + \sqrt{\tau} \bar v(n),\bar v(n)\right)}
    \end{equation}
    is Lipschitz continuous with Lipschitz constant $\sqrt{1 - \beta} M$, for some $M > 0$.
    By~\cite[Theorem 5.5]{boucheron2003concentration}, Jensen's inequality, and~\Cref{lemma:omega_stronger} we, therefore, have
    \begin{equation}\label{eq:lyapunov_discrete_stronger_final1}
        \begin{aligned}
            R_{\tau,\beta}\Vc(x,v)
            &\leq \expo{\frac{\eta^2(1-\beta)M^2}{2}} \expo{\eta\E\bigl[\sqrt{1+\omega\left(x + \sqrt{\tau} \bar v(N),\bar v(N)\right)}\bigr]}\\
            &\leq \expo{\frac{\eta^2(1-\beta)M^2}{2}} \expo{\eta(1+R_{\tau,\beta}\omega(x,v))^{1/2}}\\
            &\leq \expo{\frac{\eta^2(1-\beta)M^2}{2}} \expo{\eta\bigl(1+\omega(x,v) - \frac{\beta}{2}\|q-v\|^2 - \|q\|^2\left\{1 - \beta - \frac{L\tau}{2}\right\} + 2d(1+\frac{L\tau}{2})(1-\beta)\bigr)^{1/2}}.
        \end{aligned}
    \end{equation}
    By concavity $(1+a+s)^{1/2}\leq(1+a)^{1/2} + \frac{s}{2(1+a)^{1/2}}$ for any $a\geq 0$ and $s>-(1+a)$, thus
    \begin{equation}\label{eq:lyapunov_discrete_stronger_final2}
        \begin{aligned}
            \biggl(1+\omega(x,v)&
            -\frac{\beta}{2}\|q-v\|^2 - \|q\|^2\left\{1 - \beta - \frac{L\tau}{2}\right\} + 2d(1+\frac{L\tau}{2})(1-\beta)\biggr)^{1/2}\\
            \leq& (1+\omega(x,v))^{1/2} - \frac{\frac{\beta}{2}\|q-v\|^2 + \|q\|^2\left\{1 - \beta - \frac{L\tau}{2}\right\} - 2d(1+\frac{L\tau}{2})(1-\beta)}{2(1+\omega(x,v))^{1/2}}.
        \end{aligned}
    \end{equation}
    Note that $(q-v,q) = A(\nabla U(x), v)$ with some invertible matrix $A$ and, therefore, $\|q-v\|^2 + \|v\|^2\geq \|A^{-1}\|^{-2}\|(\nabla U(x), v))\|^2$. Using also that by assumption $\|\nabla U(x)\|^2\gtrsim \sqrt{U(x)}$ we find that
    \begin{equation}\label{eq:lyapunov_stronger1}
        \limsup_{\|(x,v)\|\rightarrow\infty}\frac{\frac{\beta}{2}\|q-v\|^2 + \|q\|^2\left\{1 - \beta - \frac{L\tau}{2}\right\} - 2d(1+\frac{L\tau}{2})(1-\beta)}{2(1+\omega(x,v))^{1/2}}>0
    \end{equation}
    Therefore, there exists $R,\rho>0$ such that for $\|(x,v)\|\geq R$
    \begin{equation}
        \begin{aligned}
            R_{\tau,\beta}\Vc(x,v)
            \leq& \expo{\frac{\eta^2(1-\beta)M^2}{2}} \expo{\eta(1+\omega(x,v))^{1/2} - \eta \rho}
            \leq \expo{\eta(1+\omega(x,v))^{1/2}} \expo{\frac{\eta(\eta (1-\beta)M^2-2\rho)}{2}}.
        \end{aligned}
    \end{equation}
    Thus, with $\eta=\rho/((1-\beta)M^2)$ to minimize the right-hand side and $\lambda \coloneqq \expo{-\frac{\rho^2}{2(1-\beta)M^2}}\in (0,1)$ we find $R_{\tau,\beta}\Vc(x,v)\leq \lambda^{1-\beta}\Vc(x,v)$.
    On the other hand, by continuity we have that $R_{\tau,\beta}\Vc(x,v)$ is bounded for $\|(x,v)\|\leq R$, concluding the proof.
\end{proof}
As mentioned above, the price we pay for the large region of allowed parameters is that the Lyapunov function degenerates as $\tau\rightarrow 0$. Indeed, one can check that the matrix $A$ such that $(q-v,v) = A(\nabla U(x),v)$ satisfies $\|A^{-1}\|^2 \gtrsim 1/\tau$. On the other hand side, the noise variance scales with $1-\beta = \Oc(\sqrt{\tau})$ under~\Cref{ass:beta} and in the case of fixed $\beta$ it does not go to zero at all which means that $R \gtrsim \tau^{-1/4}$ in the proof of~\Cref{prop:lyapunov_stronger} which tends to $\infty$ as $\tau\rightarrow 0$.
\paragraph{Lyapunov condition for small step sizes}
Similar to above we first provide the following auxiliary lemma whose proof is deferred to the supplement, \cf~\Cref{SM:lemma:ila2-recurrence,SM:lemma_lyapunov_equivalent}.
\begin{lemma}\label{lemma:ila2-recurrence}
    Let~\Cref{ass:assumptions} hold and assume that $\tau,\beta>0$ satisfy $1-L\tau-\beta^2>0$.
    Define $\omega (x, v) = U(x) + c_1\|v\|^2 + c_2\|x\|^2 + c_3\inner{x}{v}$
    with
    \begin{equation}\label{eq:discrete_lyapunov_ci}
        c_1 = \frac{1-L\tau+\beta^2}{4}, \quad
        c_2 = \frac{(1-\beta)(1-L\tau-\beta^2)}{4(1+\beta)\tau}, \quad
        c_3 = \beta\frac{1-L\tau-\beta^2}{2(1+\beta)\sqrt{\tau}}.
    \end{equation}
    Then, it holds true that
    \begin{equation}\label{eq:discrete_lyapunov}
        \begin{aligned}
            R_{\tau,\beta}\omega(x,v)\leq& \omega(x,v)
            -\frac{1-L\tau-\beta^2}{4} \|v\|^2
            -\beta\frac{1-L\tau-\beta^2}{2(1+\beta)} \inner{x}{\nabla U(x)}
            +d(1-\beta).
        \end{aligned}
    \end{equation}
    Moreover, $\omega(x,v) \asymp U(x) + \|x\|^2 + \|v\|^2$ dependently of $\tau$, $\beta$ for either $\beta$ fixed or satisfying~\Cref{ass:beta}.
\end{lemma}

    Using similar techniques as in~\Cref{prop:lyapunov_stronger} we can show the following result whose proof is deferred to the supplement~\Cref{SM:prop:lyapunov_discrete}.
    Note that now $(1-\beta)/\lambda^{1-\beta}\to 1/\log(\lambda)$ as $\beta\to 1$ rendering this Lyapunov condition stable as the step size vanishes under the assumption of \emph{kinetic momentum}~\Cref{ass:beta}.
\begin{proposition}[Drift condition for small steps]\label{prop:lyapunov_discrete}
    Let~\Cref{ass:assumptions} hold and assume that $\tau,\beta>0$ satisfy $1-L\tau-\beta^2>0$.
    Consider $\omega \colon \R^d \times \R^d \to \R$ from \Cref{lemma:ila2-recurrence}.
    Then there exist $\eta,b>0$, $\lambda\in (0,1)$ explicit in the proof and all independent of $\tau$ and $\beta$ such that for $\Vc(x,v)\coloneqq\expo{\eta\sqrt{1+\omega(x,v)}}$ it holds that
    \begin{equation}\label{eq:lyapunov-conditioon}
        R_{\tau,\beta}\Vc(x,v)\leq \lambda^{1-\beta}\Vc(x,v) + (1-\beta)b.
    \end{equation}
\end{proposition}

\begin{figure}
    \centering
    \includegraphics[]{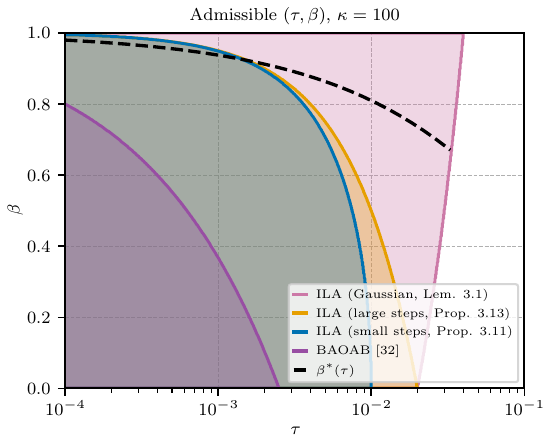}
    \caption{Depiction of the admissible parameter regions of~\eqref{eq:ila} according to~\Cref{prop:lyapunov_discrete,prop:lyapunov_stronger}, and BAOAB according to~\cite{leimkuhler2024kinetic}, as well as the optimal $\beta^*$ from~\Cref{thm_gauss}.}
    \label{fig:admissible-regions-complexity}
\end{figure}

\begin{remark}
    We want to highlight a crucial aspect at this point:~\eqref{eq:ila} does not impose lower bounds on the friction $\friction$. Indeed, the constraint $1-L\tau -\beta^2>0$ is equivalent to
    \begin{equation}
        \sqrt\tau < \frac{1-\beta^2}{\sqrt\tau L} 
        = \frac{4\friction \sqrt{\tau} + \Oc(\tau)}{\sqrt \tau L}
        = \frac{4\friction + \Oc(\sqrt\tau)}{L}
    \end{equation}
    which holds for every $\friction > 0$ once $\sqrt{\tau}$ is small enough.
    Moreover, for $\epsilon\sim \sqrt{L}$ we recover the bound $\tau\lesssim1/L$. We believe that this property crucial for acceleration.
    It has been conjectured~\cite[Section 1.1.3]{altschuler2026shifted} that among other reasons $\sqrt{\kappa}$ complexity has so far been established only for the continuous-time underdamped dynamics because acceleration requires small friction, whereas known discretizations impose a lower bound on it.
    For instance in~\cite{leimkuhler2024kinetic}, the step size condition for BAOAB reads as $\sqrt{\tau}\leq \frac{1-\expo{-2\friction \sqrt \tau}}{2\sqrt{L}}$ which, due to the inequality $1+x\leq \expo{x}$ leads to $\friction>\sqrt{L}$. Similar bounds can be found in~\cite{schuh2024convergence}. See also~\Cref{fig:admissible-regions-complexity} for a comparison of admissible parameter regions.
\end{remark}

Now, we turn to irreducibility and aperiodicity, which follow standard techniques.
We circumvent the degeneracy of the noise in~$R_{\tau,\beta}$ by considering first the two-step kernel $R_{\tau, \beta}^2$ and transferring the conclusions back to the original chain at the end.

\begin{lemma}[Existence of a continuous, strictly positive transition density]\label{lemma:density}
    The two-step transition kernel of \cref{eq:ila-2} $R_{\tau, \beta}^{2}$ admits a continuous, strictly positive density. More specifically, there exists an explicit, continuous, strictly positive function $p((x,v), (\bar x, \bar v))$ such that
    \begin{equation}
        R_{\tau, \beta}^{2}((x, v), A) = \int_A p((x,v), (\bar x, \bar v)) \dd \lambda^{2d}(\bar x, \bar v), \quad A\in \Bc(\R^{2d}).
    \end{equation}
\end{lemma}

\begin{proof}
Fix $(x,v) \in \R^{2d}$ and write $\sigma^2 = 2(1-\beta)$. Recall from \cref{eq:ila-2} that the two-step update $\Phi_{x,v} : \R^{2d} \to \R^{2d}$, $(n_1,n_2) \mapsto (\bar x,\bar v)$, is defined by
\begin{equation}\label{eq:proof-stage1}
\begin{cases}
    v_1 = \beta v - \sqrt{\tau}\, \nabla U(x) + \sigma n_1, 
    \quad 
     x_1 = x + \sqrt{\tau}\,  v_1,\\
     \bar v = \beta  v_1 - \sqrt{\tau}\, \nabla U( x_1) + \sigma n_2,
    \quad
     \bar x =  x_1 + \sqrt{\tau}\,  v_2,
\end{cases}
\end{equation}
One can easily check that $\Phi$ is continuously differentiable and $\det(D\Phi_{x,v})=(\sqrt{\tau}\sigma^2)^d$. Thus, $\Phi$ is, in fact, a diffeomorphism and $p((x,v),\cdot)$ is well-defined as the push-forward of the $2d$-dimensional standard Gaussian, denoted $\gamma^{2d}$, under $\Phi$. Specifically, we have
\begin{equation}\label{eq:trafo_rule}
    p((x,v),(\bar x,\bar v)) 
    = \left|\det (D\Phi_{x, v}^{-1}(\bar x, \bar v))\right| \gamma^{2d} (\Phi^{-1}(\bar x, \bar v))
    = \frac{\gamma^{2d} (\Phi^{-1}(\bar x, \bar v))}{(\sqrt{\tau}\sigma^2)^{d}}.
\end{equation}
Moreover, from~\eqref{eq:proof-stage1} it follows 
\begin{equation}
    \Phi^{-1}(\bar x, \bar v) = 
    \begin{pmatrix}
        \frac{(z - \sqrt{\tau} w) - x - \sqrt{\tau}\beta v + \tau \nabla U(x)}{\sqrt{\tau} \sigma}\\
        \frac{\sqrt{\tau} w - \beta ((z-\sqrt{\tau}w) - x) + \tau \nabla U(z - \sqrt{\tau} w)}{\sqrt{\tau}\sigma}
    \end{pmatrix}.
\end{equation}
which, inserted into~\eqref{eq:trafo_rule} yields the desired result. 
Continuity follows 
since $U\in C^2(\R^d)$.
\end{proof}

A strictly positive, continuous density yields a minorization condition on every compact set, and hence irreducibility and strong aperiodicity of $R_{\tau, \beta}^2$ (\cf~\cite[Section 5.4.3]{meyn1993stability}).
\begin{corollary}[Minorization condition]\label{lemma:minorization}
    For any compact set $C \in \mathcal{B}(\R^{2d})$ with non-zero Lebesgue measure there is $\eta > 0$ and a probability measure $\nu$ with $\nu(C^c)=0$ and $\nu(C)=1$, such that $R_{\tau, \beta}^2(x, A) \geq \eta \nu (A)$ for all $A \in \mathcal{B}(\R^{2d})$ and $x \in C$.
\end{corollary}
\begin{proof}
    The result follows with the choices
    \begin{equation}
        \nu(\cdot) = \frac{\lambda^{2d}(\cdot \cap C)}{\lambda^{2d}(C)},\; \eta = \bar \eta \lambda^{2d}(C)\quad \text{where }
        \bar \eta := \min_{((x,v), (\bar x, \bar v)) \in C\times C} p((x,v), (\bar x, \bar v)) > 0.
    \end{equation}
\end{proof}

With the exponential drift condition \Cref{prop:lyapunov_discrete} and the minorization condition \Cref{lemma:minorization} in hand, the standard theory now yields geometric ergodicity of \eqref{eq:ila-2}.
\begin{theorem}[Geometric ergodicity]\label{thm:discrete-ergodic}
    Let \Cref{ass:assumptions} be satisfied and let $\tau,\beta>0$ such that $\tau <\frac{2(1-\beta)}{L}$. Then $R_{\tau,\beta}$ admits a unique stationary distribution $\pi_{\tau,\beta}$ and
    there exists $C>0$, $\rho\in (0,1)$ such that for every $\mu\in \Pc(\R^{2d})$
    \begin{equation}\label{eq:disc_erg1}
        \|\mu R^k_{\tau,\beta} - \pi_{\tau, \beta}\|_{\Vc}\leq C\rho^{k} \mu(\Vc).
    \end{equation}
    Moreover, if also~\Cref{ass:beta} holds, there exists $\bar\tau>0$ such that there exist $C>0$, $\rho\in (0,1)$ such that for every $\tau<\bar \tau$ and $\mu\in \Pc(\R^{2d})$ it holds
    \begin{equation}
        \|\mu R^k_{\tau,\beta} - \pi_{\tau, \beta}\|_{\Vc}\leq C\rho^{k\sqrt{\tau}} \mu(\Vc).
    \end{equation}
\end{theorem}
\begin{proof}
    The proof follows from the standard theory of ergodicity of Markov chains, \cf~\cite[Theorem 15.0.1]{meyn2012markov}. Consider first the two-step kernel $R^{2}_{\tau,\beta}$. The existence of a strictly positive density, established in \Cref{lemma:density}, implies irreducibility.
    Moreover, by~\Cref{lemma:minorization}, every compact set $C$ with non-zero Lebesgue measure is $\nu$-small, and since $\nu(C)>0$ the chain is strongly aperiodic.
    One can easily check that the Lyapunov drift condition~\Cref{prop:lyapunov_discrete} on $R_{\tau,\beta}$ implies for $R^2_{\tau,\beta}$ the drift condition
    \begin{equation}
        R_{\tau, \beta}^2 \Vc(x,v) \leq (\lambda^2)^{1-\beta} \Vc(x, v) + (1-\beta) (1+\lambda^{1-\beta})b.
    \end{equation}
    Thus, by~\cite[Theorem 15.0.1]{meyn2012markov}, we obtain existence of a unique invariant measure $\pi_{\tau,\beta}$ and geometric ergodicity, \ie, the existence of constants $C'>0$, $\rho'\in(0,1)$ such that for any $\mu$, $\|\mu R^{2k}_{\tau, \beta} - \pi_{\tau,\beta}\|_{\Vc}\leq C' (\rho')^{k}\mu (\Vc)$. To obtain geometric ergodicity of the original chain, we consider first the \emph{shifted} two-step chain $\mu R^{2k+1}_{\tau, \beta}$. Using~\Cref{prop:lyapunov_discrete}, we get
    \begin{equation}
    \begin{aligned}
        \|\mu R^{2k+1}_{\tau, \beta} - \pi_{\tau,\beta}\|_{\Vc}
        \leq C' (\rho')^{k} \mu(R_{\tau,\beta}\Vc)
        &\leq C' (\rho')^{k} \mu(\lambda^{1-\beta}\Vc + (1-\beta)b)\\
        &\leq C' (\rho')^{k} (\mu(\Vc) + (1-\beta)b).
    \end{aligned}
    \end{equation}
    By using $R^k = (R^2)^{\frac{k}{2}}$ if $k$ is even and $R^k = R(R^2)^{\frac{k-1}{2}}$ if $k$ is odd, we obtain~\cref{eq:disc_erg1} with $\rho = \sqrt{\rho'}$ and $C = C'\vee (C'\frac{\mu(\Vc) + (1-\beta)b}{\rho \mu(\Vc)})$.

    For the second part of the assertion, note that our algorithm fits into the framework of~\cite{durmusUniformMinorizationCondition2025}. Thus, by~\cite[Theorem 3]{durmusUniformMinorizationCondition2025} we obtain a uniform minorization condition for $\tau<\bar \tau$ and some $\bar \tau>0$.
    Then, by~\cite[Theorem 5]{durmusUniformMinorizationCondition2025} we obtain the desired ergodicity.
\end{proof}
\begin{remark}
    We want to briefly compare the two regimes in the preceding proof. The first assertion applies under the explicit step size and momentum condition $1-L\tau - \beta^2>0$. The price we pay for this simple and non-restrictive condition is that the constants $C,\rho$ depend in a non-explicit way on $\tau,\beta$. However, as $\tau$ decreases, eventually we can ensure that the convergence rate behaves like $\rho^{\sqrt\tau}$ with a fixed constant, however, $\bar\tau>0$ is non-explicit.
\end{remark}


\subsubsection{Continuous Time}
We mentioned already in the introduction that the $\tau$-dependent choice of $\beta$ renders~\eqref{eq:ila} a discretization of underdamped Langevin dynamics. Specifically, under~\Cref{ass:beta} the $\tau\to 0$ limit of~\cref{eq:ila} will be
\begin{equation}\label{eq:uld}
    \begin{cases}
        \dd X_{t} = V_t \, \dd t,\\
        \dd V_t = [-2\friction V_t - \nabla U(X_t)] \, \dd t + \sqrt{4\friction} \, \wiener.
    \end{cases}
\end{equation}
In order to show that the bias of the scheme tends to zero as $\tau\rightarrow 0$ we thus have to analyze the error between~\eqref{eq:ila} and its continuous-time counterpart~\cref{eq:uld}.

While there are several results concerning ergodicity of the continuous-time dynamics, they typically either rely on coupling arguments yielding ergodicity in weaker Wasserstein metrics~\cite{cheng2018underdamped,eberle2019couplings,schuh2024convergence} or on functional inequalities yielding ergodicity in $f$-divergences~\cite{maThereAnalogNesterov2021,zhang2023improved}.
Moreover, frequently the potential is assumed to be $C^\infty$~\cite{li2026space,mattingly2002ergodicity}.
In this article, in correspondence to the discrete setting, we prove ergodicity in the $\Vc$-norm.
While the necessary Lyapunov drift condition follows from standard arguments (stated in the supplement~\Cref{SM:prop:lyapunov_cont}), in the following we provide a novel\footnote{to the best of our knowledge} proof of existence of a continuous and strictly positive density without relying on Hörmander condition. A similar result has been proven in~\cite{wu2001large}, however there the continuity is of the map $z_1\mapsto p_t(z_1,\cdot)$ as a function $\R^{2d}\rightarrow L^1(\R^{2d})$.

\begin{lemma}[Existence of a continuous, strictly positive density]
    The kinetic Langevin dynamics~\eqref{eq:uld} admit a continuous and strictly positive transition density with respect to the Lebesgue measure. More precisely, for every $t>0$ the Radon-Nikod\'ym derivative $p_t(z_1,z_2)\coloneqq \frac{\dd P_t(z_1\cdot)}{\dd \lambda^{2d}}(z_2)$
    exists and is strictly positive and jointly continuous in $(z_1,z_2)$.
\end{lemma}
\begin{proof}
    We separate the proof in multiple steps. First we derive the density of the path measure of the kinetic Langevin dynamics via Girsanov's theorem. Afterward, we use dominated convergence to show that also the time-marginal of this path measure is continuous.
    \paragraph{Step 1: Density via Girsanov}
    Omitting the potential in the kinetic Langevin equation yields the simple Ornstein-Uhlenbeck process
    \begin{equation}
        \begin{cases}
            \dd \bar X_t &= \bar V_t \dd t\\
            \dd \bar V_t &= -2\friction \bar V_t \dd t + \sqrt{2\friction}\dd \bar W_t
        \end{cases}
    \end{equation}
    initialized at $(x,v)\in\R^{2d}$. It is known that this \gls{sde} admits an explicit solution of the form
    \begin{equation}
        \bar Z_s = m_s(x,v) + L(W)_s
    \end{equation}
    with $m_s(x,v) \coloneqq (x + \friction^{-1}(1-v\expo{-\friction s}), v\expo{-\friction s})$ and
    \begin{equation}
        L(W)_s \coloneqq \int_0^s K_\sigma \dd W_\sigma
    \end{equation}
    with an appropriate deterministic function $\sigma\mapsto K_\sigma$. We will write for simplicity $L_s = L(W)_s$. In particular, the quadratic variation $\langle L\rangle_s$ reads as
    \begin{equation}
        \langle L\rangle_s=
        \begin{bmatrix}
            \frac{1}{\friction}(s - \frac{1}{\friction}(1-\expo{-2\friction s}) + \frac{1}{4\friction}(1-\expo{-4\friction s}) & \frac{1}{\friction}(1-\expo{-2\friction s}) - \frac{1}{2\friction}(1-\expo{-4\friction s})\\
            *& 1-\expo{-4\friction s}
        \end{bmatrix}.
    \end{equation}
    and we have that $\law(\bar Z_s) = \Nc(m_s(x,v),\langle L\rangle_s)$.
    Note that $m_s$ is linear, thus we may write it as $m_s(x,v)=A_s(x,v)$ with an appropriate matrix $A_s$.
    Denoting the path measure of $(\bar Z_s)_s$ started in $z = (x,v)$ as $\bar \mu^{x,v}$ and the corresponding path measure of $(X_s,V_s)_s$, the kinetic process started at $(x,v)$, as $\mu^{x,v}$, by Girsanov's theorem~\cite[Section 3.5]{karatzas1991brownian}, it follows
    \begin{equation}
        \frac{\dd \mu^{x,v}}{\dd \bar \mu^{x,v}} ((z_s)_s) = M_t ((z_s)_s) \coloneqq \expo{\frac{1}{4\friction}\int_0^t -\nabla U(x_s) \cdot(\dd v_s+2\friction v_s \dd s) - \frac{1}{8\friction}\int_0^t |\nabla U( x_s)|^2\dd s}.
    \end{equation}
    Note, however, that we are interested in a representation of the time-marginals $\mu_t^{x,v}$. Let $f$ be an arbitrary bounded measurable function. We have with $M_t = M_t((\bar Z_s)_s)$
    \begin{equation}
        \begin{aligned}
            \int f(x_t,v_t) \dd \mu_t^{x,v}(x_t,v_t)
            =\E[f(X_t,V_t)] 
            &= \E[f(\bar X_t,\bar V_t) M_t]\\
            &\hspace{-1cm}= \E[f(\bar X_t,\bar V_t) \E[M_t|(\bar X_t,\bar V_t) = (\bar X_t,\bar V_t)]]\\
            &\hspace{-1cm}= \int f(\bar x_t,\bar v_t) \E[M_t|(\bar X_t,\bar V_t) = (\bar x_t,\bar x_t)]\dd \bar \mu^{x,v}_t(\bar x_t,\bar v_t).
        \end{aligned}
    \end{equation}
    The above implies that
    \begin{equation}
        \mu_t^{x,v}(\bar x,\bar v)\dd \bar x \dd \bar v = \E[M_t|(\bar X_t,\bar V_t)=(\bar x,\bar v)]\bar\mu^{x,v}_t(\bar x,\bar v)\dd \bar x \dd \bar v,
    \end{equation}
    and since we know that $\bar\mu^{x,v}_t$ is Gaussian, whose density is smooth in $(x,v,\bar x, \bar v)$, we are left to show that 
    \begin{equation}\label{eq:time-cond-M}
        (x,v,\bar x, \bar v)\mapsto h_t(x,v,\bar x, \bar v) \coloneqq\E[M_t|(\bar X_t,\bar V_t)=(\bar x,\bar v)]
    \end{equation} 
    is continuous which we will do using dominated convergence. 
    \paragraph{Continuity of $h_t(x,v,\bar x, \bar v)$}
    We want to apply dominated convergence with respect to the underlying probability measure conditioned on the final state being $\bar Z_t = \bar z$.
    As $\bar Z_s = m_s(x,v) + L_s$ with $m_s(x,v)$ deterministic, conditioning on $\bar Z_t = \bar z$ is equivalent to conditioning on $L_t = \theta\coloneqq \bar z - m_t(x,v)$. 
    The latter is a bridge of a continuous Gaussian martingale and one can show that the distribution of $L_s|L_t = \theta$ is the same as that of an unconditioned process $L^{\theta}_t$ admitting the representation\footnote{The cited source only considers the 1D case, but we provide the arguments in the next sentence.}~\cite{gasbarra2007gaussian}
    \begin{equation}\label{eq:GP}
        \begin{aligned}
            L_s^\theta =& \langle L\rangle_s\langle L\rangle_t^{-1}\theta + \big(L_s-\langle L\rangle_s\langle L\rangle_t^{-1} L_t\big).
        \end{aligned}
    \end{equation}
    Indeed, by Gaussian conditioning the above process admits the correct marginal distributions on any finite grid $0<s_1<\dots s_n\leq t$. Since by the Kolmogorov extension theorem~\cite[Theorem 2.4.3]{tao2011introduction} this property uniquely defines a path measure, it follows that $(L^\theta_s)_s$ has the same distribution as $(L_s)_s|L_t=\theta$.\footnote{The crucial point is that not only the marginals coincide but the entire path measures.}
    Therefore, we have
    \begin{equation}
        h_t(x,v,\bar x, \bar v) 
        = \E[M_t( (\bar Z_s)_s)|\bar Z_t = \bar z]
        = \E[M_t( (m_s(x,v) + L^\theta_s)_s)].
    \end{equation}
    In other words, we can remove the conditioning by replacing $L$ with $L^\theta$. We will continue working with the unconditional version.
    
    To show continuity assume now $(x_n,v_n,\bar x_n, \bar v_n)\rightarrow (x,v,\bar x, \bar v)$ and denote 
    \begin{equation}
        K = \{(x_n,v_n,\bar x_n, \bar v_n)\;|\; n\in \N\}\cup\{(x,v,\bar x, \bar v)\}.
    \end{equation}
    To apply dominated convergence and we have to establish an integrable upper bound of $M_t$ uniformly over $K$ as well as point-wise a.e. continuity. 
    Let us rewrite the integration with respect to $\dd v_s$ in the stochastic exponential $M_t$. By the stochastic product rule
    \begin{equation}
        \dd (\nabla U(\bar X_s)\cdot \bar V_s) = \bar V_s^\top \nabla^2 U(\bar X_s) \bar V_s \dd s + \nabla U(\bar X_s)\cdot \dd V_s + \dd \langle\nabla U(\bar X_s),\bar V_s\rangle
    \end{equation}
    where we note that the crossvariation vanishes as $\bar X_s$ is defined via an \gls{ode}.
    Thus, we have 
    \begin{equation}
        \begin{aligned}
            \int_0^t \nabla U(\bar X_s)\cdot(\dd \bar V_s +2\friction \bar V_s \dd s)
            =& 2\friction\int_0^t \nabla U(\bar X_s)\cdot\bar V_s\dd s_s \\
            &\hspace{-1cm}- \int_0^t \bar V_s^\top\nabla^2 U(\bar X_s) \bar V_s \dd s + \nabla U(\bar X_t)\cdot \bar V_t - \nabla U(\bar X_0)\cdot \bar V_0.
        \end{aligned}
    \end{equation}
    It follows that
    \begin{equation}\label{eq:Mt_cont_repr}
        \begin{aligned}
            M_t
            =& \expo{\frac{1}{4\friction}\big(\nabla U(\bar X_0)\cdot \bar V_0 - \nabla U(\bar X_t)\cdot \bar V_t + \int_0^t\bar V_s^\top\nabla^2 U(\bar X_s) \bar V_s\big) - \frac{1}{2}\int_0^t\nabla U(\bar X_s)\cdot \bar V_s\dd s- \frac{1}{8\friction}\int_0^t |\nabla U( \bar X_s)|^2\dd s}.
        \end{aligned}
    \end{equation}
    Using $L$-smoothness of $U$ as well as Young's inequality multiple times yields that for some $c>0$
    \begin{equation}
        \begin{aligned}
            M_t
            \leq& \expo{(4\friction)^{-1}(\nabla U(\bar X_0)\cdot \bar V_0 - \nabla U(\bar X_t)\cdot \bar V_t}\expo{\int_0^t c(1+|\bar Z_s|^2)\dd s}\\
            \leq& \expo{(4\friction)^{-1}(\nabla U(\bar X_0)\cdot \bar V_0 - \nabla U(\bar X_t)\cdot \bar V_t)}
            \expo{tc(1+\sup_{s\leq t}|\bar Z_s|^2)}
        \end{aligned}
    \end{equation}
    for some $c>0$. Let us write
    $L^\theta$ as $L_t^\theta = a_t(\theta) + \xi_t$ with $a_t(\theta)$ deterministic and $\xi_t$ a zero mean Gaussian process independent of $\theta$ with a.s. continuous paths. Thus, we have
    \begin{equation}\label{eq:upperbound_MT}
        \begin{aligned}
            \E[\sup_{(z,\bar z)\in K} M_t( (m_s(x,v) + L^\theta_s)_s)]
            &\\
            &\hspace{-4cm}\leq\expo{(4\friction)^{-1}(\nabla U(x)\cdot v-\nabla U(\bar x)\cdot \bar v)}\times \E[\sup_{(z,\bar z)\in K}\expo{tc(1+\sup_{s\leq t}|m_s(x,v) + a_s(\theta) + \xi_t|^2)}]\\
            &\hspace{-4cm}\leq \expo{ct}\expo{c(\nabla U(x)\cdot v-\nabla U(\bar x)\cdot \bar v + ct\sup_{s\leq t}(|m_s(x,v)|^2 + |a_s(\theta)|^2)}\E[\expo{tc \sup_{s\leq t}|\xi_t|^2)}]
        \end{aligned}
    \end{equation}
    where $c>0$ is an absolute constant.
    Of course the Gaussian process $\xi_s$ admits Gaussian tails. More precisely, by Fernique's theorem~\cite[Remark 4.15]{hairer2009introduction}, $\E[\expo{tc \sup_{s\leq t}|\xi_t|^2)}]<\infty$ as long as $ct<\frac{1}{2\sigma^2}$ where
    \begin{equation}
        \sigma^2 = \sup_{s\leq t}\E[|\xi_s|^2].
    \end{equation}
    One can easily verify from~\eqref{eq:GP} that $\sigma^2 \rightarrow 0$ as $t\rightarrow 0$.
    Since $(\xi_s)_s$ is, moreover, independent of $(x,v,\bar x, \bar v)$, we can choose $t$ sufficiently small so that $ct<\frac{1}{2\sigma^2}$ also independent of $(x,v,\bar x, \bar v)$. As a consequence, the expected value on the right-hand side in~\eqref{eq:upperbound_MT} is bounded independent of $(x,v,\bar x, \bar v)$ and the other terms are bounded for $(x,v,\bar x, \bar v)\in K$ due to compactness and continuity. Thus, we have shown the uniform upper bound needed to apply dominated convergence. 
    Pointwise-a.e. convergence follows directly from a.e. continuity of the involved functions.
    We may therefore apply dominated convergence which shows that $\mu^{x,v}_t(\bar x, \bar v)$ admits a density which is continuous jointly in $(x,v,\bar x, \bar v)$ as long as $t$ is sufficiently small. 
    \paragraph{Extending the result to arbitrary $t>0$}    
    In order to extend the result to arbitrary $t$, fix $t_0$ small enough so that the result is valid for $t=t_0$. We will show inductively that $\mu^{x,v}_{kt_0}(\bar x, \bar v)$ is continuous jointly in $(x,v,\bar x, \bar v)$ from which the result can then be derived by expressing arbitrary $t>0$ as $t=kt_0+\delta$, $\delta \in [0,t_0)$. More specifically, we show that for each $k$ we have
    \begin{enumerate}
        \item $\mu^{x,v}_{kt_0}(\bar x, \bar v)$ is continuous jointly in $(x,v,\bar x, \bar v)$ and
        \item For any compact set $K\subset \R^{2d}$ there exists $C(t_0,K)>0$ such that
        \begin{equation}
            \begin{aligned}
                \sup_{(x,v)\in \R^{2d}, (\bar x, \bar v)\in K}\mu^{x,v}_{kt_0}(\bar x, \bar v) \leq C(t_0,K)\\
            \end{aligned}
        \end{equation}
    \end{enumerate}
    In the first part of the proof we have established the first assertion for $k=1$.
    To prove the second one, note that the covariance of the distribution $\bar \mu^{x,v}_s$, $\langle L\rangle_s$ converges to zero as $s\rightarrow 0$ independently of $(x,v)$.
    Therefore, by potentially decreasing $t_0$ we can ensure that for any compact $K$ there exists a constant $C(t_0,K)>0$ such that
    \begin{equation}\label{eq:sup_mu}
        \begin{aligned}
            \sup_{(x, v)\in\R^{2d}, (\bar x, \bar v)\in K}\mu^{x, v}_{t_0}(\bar x, \bar v)
            =& \sup_{(x, v)\in\R^{2d}, (\bar x, \bar v)\in K}\bar \mu^{x, v}_{t_0}(\bar x, \bar v)h_{t_0}(x, v,\bar x, \bar v)\\
            \leq& \sup_{(x, v)\in\R^{2d}, (\bar x, \bar v)\in K}\bar \mu^{x, v}_{t_0}(\bar x, \bar v)\expo{ct_0}\\
            &\times\expo{c(\nabla U(x)\cdot v - \nabla U(\bar x)\cdot \bar v+ ct_0\sup_{s\leq t_0}(|m_s(x,v)|^2 + |a_s(\theta)|^2))}\E[\expo{ct_0 \sup_{s\leq t_0}|\xi_t|^2}]\\
            \leq& C(t_0,K)<\infty.
        \end{aligned}
    \end{equation}
    Indeed, this follows from the fact that for any $\delta>0$ we can find $t_0$ small enough such that 
    \begin{equation}\label{eq:sup_mu2}
        \bar \mu_{t_0}^{x,v}(\bar x, \bar v)
        \leq \frac{1}{(2\pi\det(\Sigma_{t_0})^{1/2}}\expo{-\frac{\|(\bar x, \bar v)-m_{t_0}(x,v)\|^2}{2\delta}}
        \leq \frac{1}{(2\pi\det(\Sigma_{t_0}))^{1/2}}\expo{-\frac{(\|m_{t_0}(x,v)\| - \|(\bar x, \bar v)\|)^2}{2\delta}}.
    \end{equation}
    Since also $t_0\sup_{s\leq t_0}(|m_s(x,v)|^2 + |a_s(\theta)|^2))$ remains bounded as $t_0\rightarrow 0$, we can choose $\delta>0$ sufficiently small so that the decay of $\bar \mu_{t_0}^{x,v}(\bar x, \bar v)$ compensates any quadratic growth in~\eqref{eq:sup_mu} leading to boundedness. 

    For the induction step first note that we can write
    \begin{equation}\label{eq:Kolmogorov}
        \begin{aligned}
            \mu^{x,v}_{(k+1)t_0}(\bar x, \bar v) =& \int \mu^{x,v}_{t_0}(\dd x',\dd v') \mu^{x', v'}_{kt_0}(\bar x, \bar v)\\
            =& \int \bar \mu^{x,v}_{t_0}(\dd x',\dd v') h_{t_0}(x,v,x', v')\mu^{x', v'}_{kt_0}(\bar x, \bar v).
        \end{aligned}
    \end{equation}
    Let $(x_n,v_n,\bar x_n, \bar v_n)\rightarrow (x,v,\bar x, \bar v)$. The integrand in~\eqref{eq:Kolmogorov} is continuous for a.e. $(x', v')$ by the induction hypothesis so that we are only left to show again that it admits an integrable upper bound uniformly in $n$ to use dominated convergence. By the induction hypothesis we have $\sup_{n} \mu^{x', v'}_{k t_0}(x'_n,v'_n)\leq C <\infty$ uniformly over $(x', v')$. Moreover, the same arguments as used in~\eqref{eq:sup_mu}, \eqref{eq:sup_mu2} yield that $\sup_n\bar \mu^{x_n,v_n}_{t_0}(\dd x',\dd v') h_{t_0}(x_n,v_n,x', v')$ is integrable over $(x',v')$. Thus, dominated convergence applies and we obtain continuity. The second assertion easily follows as for any compact $K$ and by the induction hypothesis
    \begin{equation}
        \begin{aligned}
            \sup_{(x,v)\in \R^{2d}, (\bar x, \bar v)\in K}\mu^{x,v}_{(k+1)t_0}(\bar x, \bar v) 
            =& \sup_{(x,v)\in \R^{2d}, (\bar x, \bar v)\in K}\int \mu^{x,v}_{t_0}(\dd x',\dd v') \mu^{x', v'}_{kt_0}(\bar x, \bar v)\\
            \leq& \sup_{(x,v)\in \R^{2d}, (\bar x, \bar v)\in K}\int \mu^{x,v}_{t_0}(\dd x',\dd v') C(t_0,K)\\
            =& C(t_0,K)\\
        \end{aligned}
    \end{equation}    
    concluding the proof.
\end{proof}

\begin{theorem}[Exponential ergodicity]\label{thm:cont_ergodic}
    There exist $C>0$ and $\kappa\in (0,1)$ such that for any initial distribution $\mu$, $\|\mu P_t - \pi \|_{\Vc}\leq C \mu(\Vc)\kappa^t$
    and, thus, since $\Vc\geq 1$ also $\|\mu P_t - \pi\|_{\TV}\leq C \mu(\Vc)\kappa^t$ and for every $p\geq1$ there exists $C_p>0$ such that $\Wc_p(\mu P_t,\pi)\leq C_p \mu(\Vc)^{1/p}\kappa^t$.
\end{theorem}
\begin{proof}
    The strictly positive and continuous density implies as in the discrete case that for every skeleton chain (\ie, every chain with transition kernel $P_\delta$ for a $\delta>0$), every compact set with non-zero Lebesgue measure is small. Then the Lyapunov condition implies the desired result by~\cite[Theorem 6.1]{meyn1993stability}.
    The convergence in Wasserstein then follows directly by~\cite[Theorem 6.15]{villani2008optimal} and since for every $p\geq 1$ there exists $C_p$ such that $\|(x,v)\|^p\leq C_p\Vc(x,v)$.
\end{proof}

\subsubsection{Bounding the discretization error}

We conclude the theoretical part by showing that also in the setting of general potentials $U$ we have $\pi_{\tau,\beta}\rightarrow \pi$ as $\tau\rightarrow 0$. In the following, we switch to the natural time-scale for underdamped Langevin, that is, we denote $\step\coloneqq \sqrt{\tau}$. Recall that $P_t$ denotes the Markov semigroup of~\eqref{eq:uld}. The proof of the following relies on standard techniques and is, thus, deferred to the supplement, see~\Cref{SM:thm:disc_SM}
\begin{proposition}[Discretization error]\label{thm:disc}
    Let~\Cref{ass:assumptions,ass:beta} hold true and denote $\step\coloneqq\sqrt{\tau}$. Then there exists $\rho, c >0$ such that for any initial distribution $\mu\in\Pc_2(\R^{2d})$ it holds true
    \begin{equation}
        \Wc_2^2(\mu R_{\tau,\beta}^k,\mu P_{k\step})\leq k \disca \step^2\expo{c k\step},\quad k\geq 0.
    \end{equation}
\end{proposition}
\begin{theorem}
    It holds $\Wc_2(\pi_{\tau\beta},\pi)\rightarrow 0$ as $\tau\rightarrow 0$.
\end{theorem}
\begin{proof}
    By the triangle inequality and~\Cref{thm:cont_ergodic}, it holds
    \begin{equation}
        \begin{aligned}
            \Wc_2(\pi_{\tau,\beta} R^k_{\tau,\beta},\pi) =& \Wc_2(\pi_{\tau,\beta} R^k_{\tau,\beta},\pi)\\
            \leq& \Wc_2(\pi_{\tau,\beta} R^{k_1}_{\tau,\beta} R^{k-k_1}_{\tau,\beta}, \pi_{\tau,\beta} R^{k_1}_{\tau,\beta} P_{(k-k_1)\step}) + \Wc_2(\pi_{\tau,\beta} R^{k_1}_{\tau,\beta} P_{(k-k_1)\step},\pi)\\
            \leq& \sqrt{\rho (k-k_1)\step^2\expo{c(k-k_1)\step}} + C_2 \pi_{\tau,\beta}(\Vc)^{1/2}\kappa^{(k-k_1)\step}.
        \end{aligned}
    \end{equation}
    By~\Cref{prop:lyapunov_discrete}, $\pi_{\tau,\beta}(\Vc)$ is bounded by a finite constant uniformly in $\tau$.
    The result follows upon choosing $T\coloneqq (k-k_1)\step$ sufficiently large to make the second term small and afterward $\step$ sufficiently small to make the first term small.
\end{proof}

\begin{remark}
    Using similar techniques as in~\cite[Theorem 7.3]{mattingly2002ergodicity} one can also derive existence of $\delta>0$ such that for any test function $g:\R^{2d}\rightarrow \R$, with appropriate local Lipschitz continuity $|\pi_{\tau,\beta}(g) - \pi(g)|\lesssim \tau^\delta$.
\end{remark}

\section{Numerical Experiments}\label{sec:numerical}
In this section, we showcase the applicability of our algorithm~\eqref{eq:ila} in two real-world use cases.
Specifically, we employ~\eqref{eq:ila} on a high-dimensional Bayesian image denoising task in~\Cref{ssec:ex-tv-denoising} and perform maximum likelihood training for an \gls{ebm} in the context of molecular structure generation in~\Cref{ssec:ex-mol-gen}.

\subsection{TV Denoising}\label{ssec:ex-tv-denoising}
\begin{figure}
    \centering
    \includegraphics[scale=0.8]{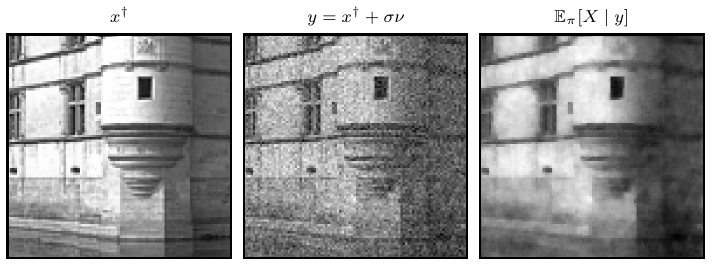}
    \caption{We see the clean ground truth signal (left), the noisy observation serving as conditioning (middle), and the \gls{mmse} estimate of the ground truth signal given the noisy observation approximated by the sample mean over the ground truth samples drawn with~\cite{kuric_gaussian_2025}.}
    \label{fig:ex-tv-denoising-gt-noisy-mmse}
\end{figure}

We consider the inverse problem of recovering a $d = d_1 \times d_2 = 96 \times 96$ dimensional grayscale image from its noisy observation (see~\Cref{fig:ex-tv-denoising-gt-noisy-mmse}).
To tackle this problem we use the ROF model~\cite{rudin_nonlinear_1992} $U(x) = \frac{1}{2\sigma^2}\|x-y\|^2 + \lambda \TV(x)$ with an observation $y\in\R^d$, a noise level $\sigma>0$, a regularization parameter $\lambda>0$ and the anisotropic total variation functional 
\begin{equation}
    \TV(x) = \sum_{i=1}^{d_1} \sum_{j=1}^{d_2} \abs{(k_h \ast P x)_{i,j}} + \abs{(k_v \ast Px)_{i,j}},
\end{equation}
where $P \colon \R^{d} \to \R^{(d_1 + 2 ) \times (d_2 + 2)}$ denotes a circular padding operator and the convolution kernels $k_h = [0,-1,1]$ and $k_v = k_h^\top$
extract horizontal and vertical forward finite differences.
The potential $U$ is the negative log-posterior under the measurement model $y = x + \sigma \nu$, $\nu \sim \Nc(0, \id_d)$ when assuming a prior distribution $x\sim p(x)\propto \expo{-\lambda \TV(x)}$\footnote{Note that this is, in fact, only a proper (=integrates to one) probability distribution on a subspace of $\R^d$. The posterior, however, is proper due to the Gaussian term.}.
Since the \gls{tv} functional is non-smooth, the potential in its current form would violate the theoretical assumptions.
Thus, we resort to approximating the absolute value with the Huber function, which gives us a strongly convex, $L$-smooth potential $U$ meeting \Cref{ass:assumptions}.

\begin{figure}
    \centering
    \includegraphics[]{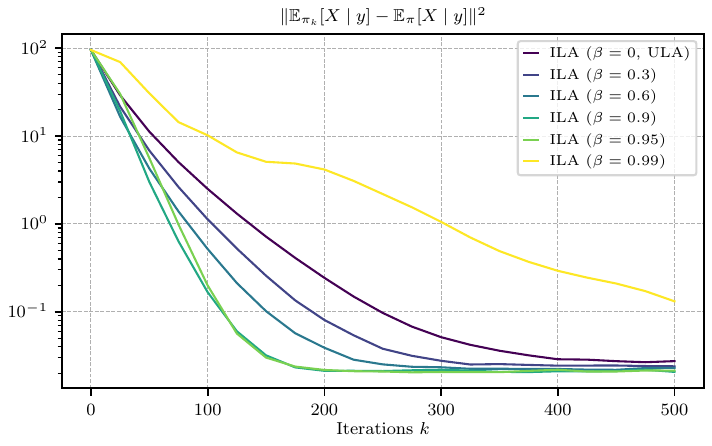}
    \caption{Squared error between estimated posterior mean and the ground truth \gls{mmse} (depicted in \Cref{fig:ex-tv-denoising-gt-noisy-mmse} on the right) over the iterations.
    Momentum eases convergence across a wide range of values for $\beta$.}
    \label{fig:ex-tv-denoising-mmse-diff}
\end{figure}

Our goal now is to test the effect of momentum on the sampling speed.
Thus, we choose $\beta \in \{0, 0.3, 0.6, 0.9, 0.95, 0.99\}$ and set $\tau(\beta) = 0.85 \cdot (1-\beta^2)/L$ in accordance to \Cref{thm:discrete-ergodic} for each instance of \eqref{eq:ila}.
For each run, we produce $N=\num{1000}$ reconstructions by simulating the dynamics initialized in the noisy image for $K=500$ iterations.
To assess the speed of convergence we plot the deviation of the \gls{mmse} estimate produced by \eqref{eq:ila} to the ground truth \gls{mmse} in \Cref{fig:ex-tv-denoising-mmse-diff}.
It is evident, that adding momentum speeds up convergence.
The ground truth samples are produced with the \gls{glm}~\cite{kuric_gaussian_2025} with $\sigma = 0.1$ and $\lambda = 12.58714$ (the optimal regime reported in~\cite{kuric_gaussian_2025}).
The final \gls{mmse} estimates produced by the \eqref{eq:ila} instances can be found in \Cref{SM:sec:tv-denoising}.

\begin{figure}
    \centering
    \includegraphics[]{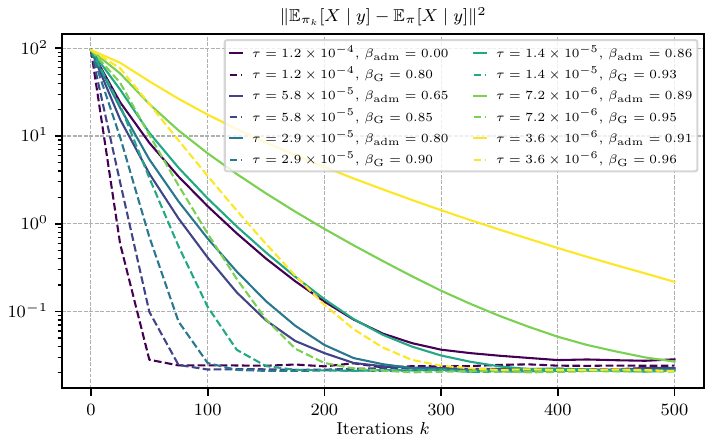}
    \caption{Squared error between estimated posterior mean and the ground truth \gls{mmse} (depicted in \Cref{fig:ex-tv-denoising-gt-noisy-mmse} on the right) over the iterations. For each value of $\tau$ $\beta_{\mathrm{adm}}$ is a $\beta$ admissible according to~\Cref{prop:lyapunov_discrete} and $\beta_{\mathrm{G}} = \beta^\ast(\tau)$ the optimal choice according to~\Cref{thm_gauss}.
    The Gaussian-optimal choice yields acceleration even though the target is not Gaussian.}
    \label{fig:ex-tv-denoising-mmse-diff-tau-sweep}
\end{figure}

To test whether our analysis in the general case (\cf~\Cref{sec:non_gauss}) leaves room for improvement, we also perform a sweep over step sizes $\tau_i = 1 / (2^i L), i = 0, \dots, 5$ and compare the resulting momentum from the general analysis $\beta_{\text{adm}} < \sqrt{1 - L\tau}$ with the rate-optimal one from the Gaussian case $\beta_{\text{G}} = (1 - \sqrt{m \tau})^2$, where in our case $m = 1 / \sigma^2$.
From \Cref{fig:ex-tv-denoising-mmse-diff-tau-sweep}, it is evident that $\beta_G$ consistently outperforms $\beta_{\text{adm}}$ indicating that the convergence results of the Gaussian case should be true also in the strongly convex setting.

\subsection{Molecular Structure Generation}\label{ssec:ex-mol-gen}
\begin{figure}
    \centering
    \includegraphics[width=0.49\linewidth]{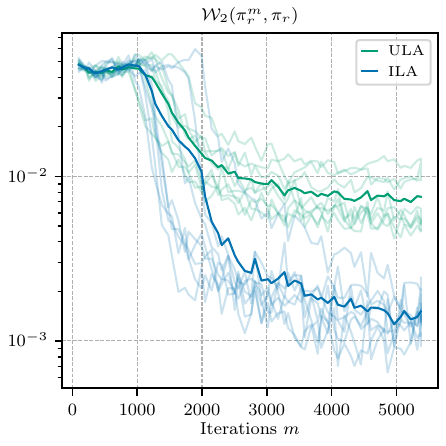}
    \hfill
    \includegraphics[width=0.49\linewidth]{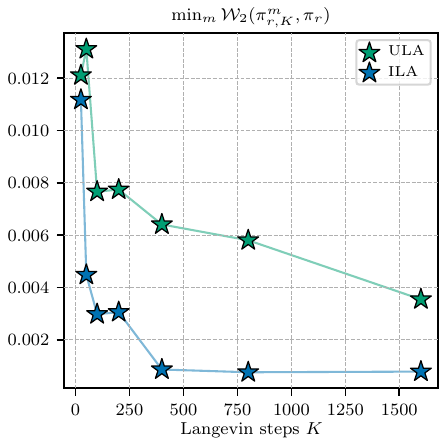}
    \caption{\emph{Left:} Convergence rate of Wasserstein-2 distance between bond length distributions for \gls{ula} and~\eqref{eq:ila} over training iterations $m$ and a constant number of $K=400$ Langevin steps per iteration.
    We show the results for multiple training runs, the mean performance is depicted using solid lines.
    \emph{Right:} Wasserstein-2 distance between bond length distributions when varying the number of Langevin steps $K$ performed in each iteration of the training. We pick the best-performing model checkpoint over $m$ iterations.}
    \label{fig:ex-mol-gen}
\end{figure}
Here, we apply~\eqref{eq:ila} to the challenging task of maximum likelihood learning of a Lipschitz-constrained \gls{ebm}. That is, in this experiment we do not evaluate the \emph{isolated} sampling performance, but instead how well a sampler performs when used within a training task.
Specifically, we address the problem of \emph{molecular structure generation}~\cite{baillif_deep_2023}, a fundamental problem in computational chemistry and structural biology.
The goal is to predict energetically favorable 3D configurations for a given set of atoms.
\Glspl{ebm} allow for a natural treatment of this problem by directly parameterizing the molecular \gls{pes}.
Such energy landscapes are known to be highly multi-modal, with many shallow local minima, and physically valid conformations occupy only a small fraction of the configuration space~\cite{axelrod_geom_2022}.
As \gls{ebm} training via contrastive divergence~\cite{hinton_training_2002} crucially depends on the quality of samples drawn from the model, this task serves as a challenging benchmark for the sampler's ability to traverse complex, multi-modal potential functions.
In the following, we demonstrate the ability of~\eqref{eq:ila} to improve mixing and convergence in such multi-modal settings, and benchmark it against the predominant \gls{ula} algorithm.

We are interested in learning the parameters $\psi$ of a conditional Gibbs distribution $\pi(x \mid a) \propto \exp\bigl(-E_\psi(x \mid a)\bigr)$,
where $x = \{x_i\}_{i=1}^V \in \mathbb{R}^{V \times 3}$ are the 3D coordinates of the $V$ atoms, and $a = \{a_i\}_{i=1}^V \in \mathcal{V}^V$ their atomic types drawn from the set $\mathcal{V}$.
The potential function that we model using a neural network with parameters $\psi$ is denoted via $E_\psi(x \mid a) \in \mathbb{R}$. 

We parameterize the energy function $E_\psi$ with a three-layer \gls{egnn}~\cite{satorras_en_2021} with embedding dimension $\num{128}$.
The final energy is obtained by summing the per-atom outputs of the last layer.
To better leverage our theoretical results regarding valid momentum parameters, we apply spectral normalization~\cite{miyato_spectral_2018} to all linear layers of the \gls{egnn}, effectively upper-bounding the network's Lipschitz constant with $L \leq 1$.
Conveniently, this addition is also known to increase training stability~\cite{du_implicit_2020}.
The training data consists of a subset of the QM9 dataset~\cite{ramakrishnan_quantum_2014}, a popular dataset for the molecule generation task.
To reduce the computational workload, we restricted the dataset to molecules of size $V=19$.
The resulting dataset contains $\num{18146}$ molecular structures, which we split into $90\%$ training and $10\%$ test data.
Bond types are not provided during training, meaning the network only receives atom types and inter-atomic distances as input.
As the primary evaluation criterion, we calculate how well the model captures the chemically plausible bond length distribution.
Therefore, we compute the pairwise distance distribution of atoms $\pi_r^m$ of our model at iteration $m$ and compare it to the data distribution.
As we are mainly interested in accurate modeling of bond lengths, we restrict the considered pairwise distances to a range from \SIrange{1}{2}{\angstrom} and compare to the ground-truth distribution $\pi_r$ using the Wasserstein-2 distance $\mathcal W_2(\pi_{r}^m, \pi_{r})$.
More details about the training procedure of the \gls{ebm} can be found in \Cref{SM:app:ebm_training}.

In \Cref{fig:ex-mol-gen}, the left plot reports the Wasserstein-2 distance between the bond length distributions of the test set and model samples during training, using a fixed budget of $K=400$ Langevin steps.
While both samplers initially improve, \gls{ula} rapidly saturates at a higher distance value, indicating that its chains fail to mix once the energy landscape becomes increasingly rugged through contrastive divergence updates.
In contrast,~\eqref{eq:ila} maintains stable exploration, traversing shallow minima and continuing to provide diverse and informative samples.
This allows the learned energy function to improve throughout training, yielding substantially lower $\mathcal {W} _ 2$-distances at the same computational cost.

We next analyze how the computational budget influences model quality.
The right plot in~\Cref{fig:ex-mol-gen} shows the final model quality over the available computational budget, where each point corresponds to the best checkpoint of an independent run.
Increasing the number of Langevin steps improves both methods, yet ILA exhibits a much steeper initial gain and reaches its asymptotic accuracy after only a fraction of the iterations required by \gls{ula}.
With only $\num{100}$ inner-loop steps,~\eqref{eq:ila} already matches the final performance of \gls{ula}, corresponding to a roughly $\num{16}$-fold reduction in computational cost, and converges fully at around $\num{400}$ steps.
The addition of momentum substantially enhances sampling efficiency, particularly in curvature-constrained energy landscapes, such as in our case, where spectral normalization bounds the gradient norm and thus limits local curvature.
This reduces the need for long inner-loop simulations and effectively mitigates one of the primary computational bottlenecks in maximum-likelihood training of energy-based models.

\section{Related Work}\label{sec:related}
\paragraph{Overdamped Langevin}
There is a vast body of work regarding overdamped Langevin dynamics~\eqref{eq:OD} for sampling using with proof covering Wasserstein coupling arguments~\cite{eberle2016reflection}, Lyapunov drift conditions~\cite{durmus2017nonasymptotic,durmus2019high,robertsExponentialConvergenceLangevin1996}, functional inequalities~\cite{chewi2025analysis,vempala2019rapid} or optimization in Wasserstein space~\cite{durmus2019analysis}. For a comprehensive overview over different proof strategies, see also~\cite{chewi2023log}.

\paragraph{Underdamped Langevin}
Due to its accelerating properties there has been substantial interest in underdamped Langevin dynamics for sampling in recent years.
However, proving $\sqrt{\kappa}$ complexity has proven difficult. 
Similarly to the overdamped case, different strategies have been proposed and the literature is too exhaustive to give a comprehensive overview. Due to the degeneracy of the Brownian motion, the proofs are often more involved than in the case of overdamped Langevin dynamics which, in particular, incentivised the development of hypocoercivity~\cite{villani2006hypocoercivity}.
An early result in Wasserstein-2 using synchronous coupling and leading to $\Oc(\sqrt{d}\kappa^2/\delta)$ complexity is provided in~\cite{cheng2018underdamped}. Using similar techniques, the results were improved in~\cite{dalalyan2020sampling} where, however, for large $\kappa$ one still obtains complexity scaling with $\kappa^2$.
In~\cite{eberle2019couplings,schuh2024convergence} reflection coupling together with specifically designed distance functions are used to show ergodicity in Wasserstein metrics. In the latter work, the authors obtain in the Gaussian similarly scaling of $\Oc(\kappa L)$.
\cite{leimkuhler2024kinetic} give a comprehensive treatment of contraction results of various discretizations of underdamped Langevin in Wasserstein two using synchronous coupling. For BAOAB they obtain $\Oc(\kappa\sqrt{L})$ scaling
For a treatment of splitting methods using functional analytic tools, see~\cite{leimkuhler2016computation}.
\cite{durmusUniformMinorizationCondition2025} provide a proof relying on classical Lyapunov analysis similar to the present article in the non-Gaussian setting.
For analyses in Kullback-Leibler and Chi-squared distance relying on functional inequalities, see~\cite{maThereAnalogNesterov2021,zhang2023improved}

Recently, efforts have been made to close the gap to the desired $\sqrt{\kappa}$ complexity~\cite{zhang2023improved,altschuler2026shifted}. In the continuous time setting, $\Oc(\sqrt{\kappa})$ has been reached via space-time log-Sobolev constants~\cite{cao2023explicit}. Most notably, in~\cite{altschuler2026shifted} for the first time sublinear in $\kappa$ complexity, namely $\Oc(\kappa^{5/6})$, has been achieved for strongly convex and $L$-smooth potentials. As mentioned in the introduction~\cite{bou2026provable} provides a negative result for a range of splitting methods.

\section{Conclusion}
In this article, we introduced~\eqref{eq:ila}, an algorithm for sampling from Gibbs distributions that mirrors Polyak's heavy-ball method from optimization.
We showed that it attains $\sqrt{\kappa}$ complexity in the Gaussian setting, reaching the ballistic rate known from optimization.
In the general setting we establish geometric ergodicity of the discrete scheme and exponential ergodicity of its continuous-time counterpart, together with convergence of the discretization bias to zero as the step size decreases -- all under considerably simpler
parameter restrictions than previously available, notably without any lower bound on the
friction coefficient.
The practical relevance of~\eqref{eq:ila} is confirmed by extensive numerical experiments, confirming the $\sqrt{\kappa}$ scaling in the Gaussian case and demonstrate the efficiency of the method in complex machine learning scenarios.

\paragraph{Limitations and future work}
The greatest limitation and, thus, most interesting aspect to consider in future work is an improvement of the analysis in the non-Gaussian setting.
Additionally, we would like to analyze a similar scheme that resembles Nesterov's gradient descent as we believe this might be a key to closing the gap to $\sqrt\kappa$ scaling.

\bibliographystyle{siamplain}
\bibliography{references}

\includepdf[fitpaper=true, pages=-]{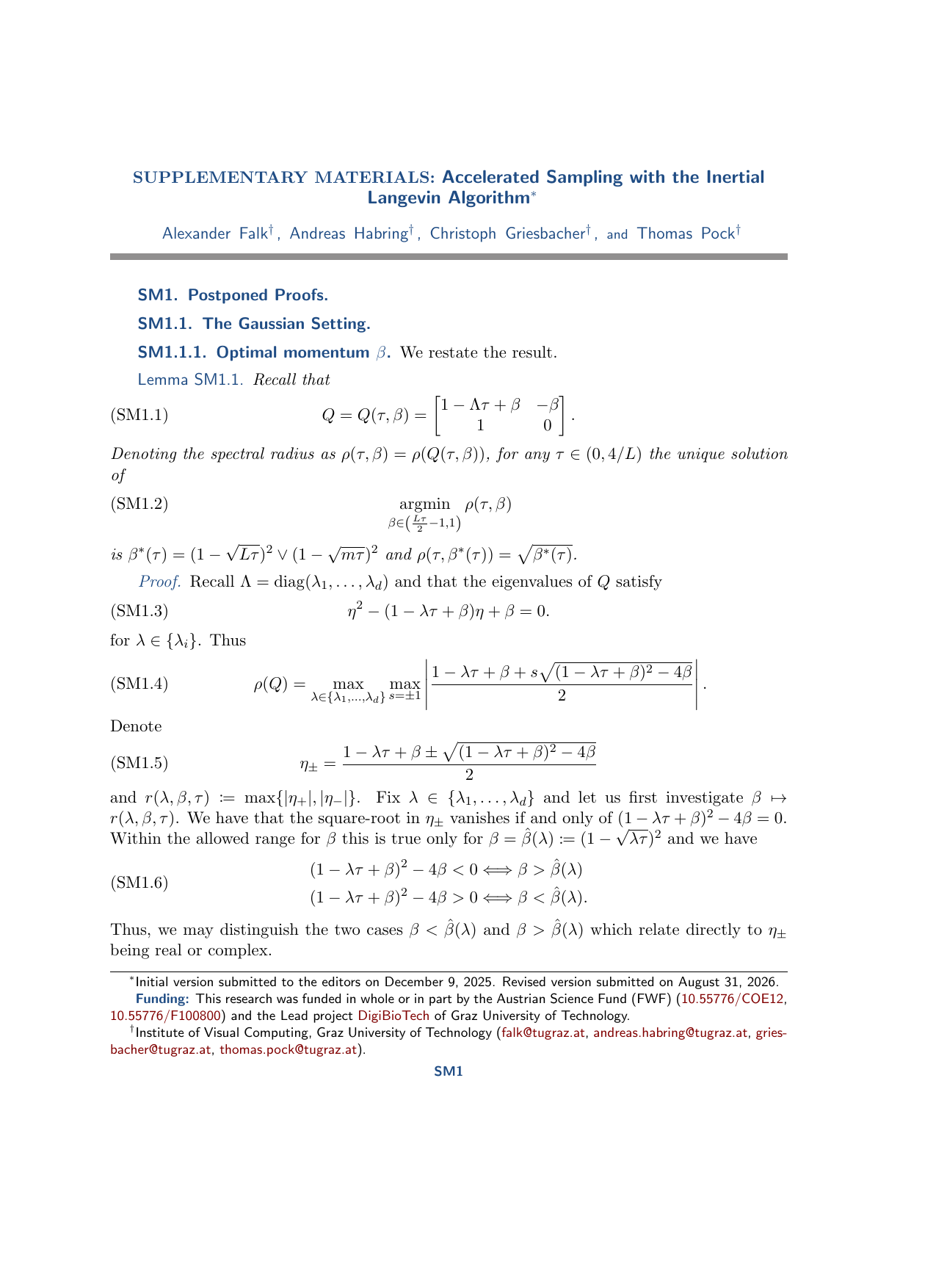}
\makeatletter\@input{xx2.tex}\makeatother
\end{document}


\maketitle
\section{Postponed Proofs}
\subsection{The Gaussian Setting}
\subsubsection{Optimal momentum $\beta$}
We restate the result.
\begin{lemma}\label{lemma:supp:optimal_beta}
        Recall that
        \begin{equation}
            Q=Q(\tau,\beta)=\begin{bmatrix}
            1-\Lambda\tau + \beta&-\beta\\
            1&0
        \end{bmatrix}.
        \end{equation}
        Denoting the spectral radius as $\rho(\tau, \beta) = \rho (Q(\tau, \beta))$, for any $\tau\in (0,4/L)$ the unique solution of
        \begin{equation}
            \argmin_{\beta\in \left(\frac{L\tau}{2}-1,1\right)}\rho(\tau, \beta)
        \end{equation}
        is $\beta^*(\tau) = (1-\sqrt{L\tau})^2\vee (1-\sqrt{m\tau})^2$ and $\rho(\tau,\beta^*(\tau)) = \sqrt{\beta^*(\tau)}$.
    \end{lemma}
\begin{proof}
    Recall $\Lambda = \diag(\lambda_1,\dots,\lambda_d)$ and that the eigenvalues of $Q$ satisfy
    \begin{equation}\label{eq:heavy_ball_convergence}
        \eta^2-(1-\lambda\tau + \beta)\eta + \beta  = 0.
    \end{equation}
    for $\lambda \in \{\lambda_i\}$. Thus
    \begin{equation}
        \begin{aligned}
            \rho(Q) 
            = \max_{\lambda\in \{\lambda_1,\dots,\lambda_d\}}\max_{s=\pm 1}\left|\frac{1-\lambda\tau + \beta + s\sqrt{(1-\lambda\tau + \beta)^2 - 4\beta}}{2}\right|.
        \end{aligned}
    \end{equation}
    Denote 
    \begin{equation}\label{eq:eta}
        \eta_{\pm} = \frac{1-\lambda\tau + \beta \pm \sqrt{(1-\lambda\tau + \beta)^2 - 4\beta}}{2}
    \end{equation}
    and $r(\lambda,\beta,\tau) \coloneqq \max\{|\eta_+|,|\eta_-|\}$.
    Fix $\lambda\in \{\lambda_1,\dots,\lambda_d\}$ and let us first investigate $\beta\mapsto r(\lambda,\beta,\tau)$. 
    We have that the square-root in $\eta_\pm$ vanishes if and only of $(1-\lambda\tau+\beta)^2-4\beta=0$. Within the allowed range for $\beta$ this is true only for $\beta = \hat\beta(\lambda)\coloneqq (1-\sqrt{\lambda\tau})^2$ and we have
    \begin{equation}
        \begin{aligned}
            (1-\lambda\tau + \beta)^2 - 4\beta<0 &\Longleftrightarrow \beta>\hat\beta(\lambda)\\
            (1-\lambda\tau + \beta)^2 - 4\beta>0 &\Longleftrightarrow \beta<\hat\beta(\lambda).
        \end{aligned}
    \end{equation}
    Thus, we may distinguish the two cases $\beta<\hat\beta(\lambda)$ and $\beta>\hat\beta(\lambda)$ which relate directly to $\eta_\pm$ being real or complex.
    \paragraph{Case 1: $\beta>\hat\beta(\lambda)$}
    Since $\eta_\pm$ is complex we have
    \begin{equation}
        \begin{aligned}
            r(\lambda,\beta) = |\eta_+| = |\eta_-| = \sqrt{\beta}
        \end{aligned}
    \end{equation}
    In particular, note that $r(\lambda,\beta)$ is independent of $\lambda$ and strictly increasing in $\beta$.
    \paragraph{Case 2: $\beta<\hat\beta(\lambda)$}
    In this case, $\eta_\pm\in \R$.
    Recall that $\eta_{\pm}$ solves~\eqref{eq:heavy_ball_convergence}. Therefore, considering $\eta_\pm = \eta_\pm(\beta)$ as functions of $\beta$, which are continuously differentiable for $\beta<\hat\beta$, we obtain
    \begin{equation}
        0 = \frac{\dd}{\dd\beta }(\eta_\pm^2-(1-\lambda\tau + \beta)\eta_\pm + \beta) = 2\eta_\pm\eta_\pm' -\eta_\pm - (1-\lambda\tau+\beta)\eta_\pm'+1
    \end{equation}
    and therefore, $\eta_\pm' = \frac{\eta_\pm-1}{2\eta_\pm-(1-\lambda\tau+\beta)}$, that is,
    \begin{equation}\label{eq:heavy_ball_derivatives}
        \eta_\pm' = \frac{\eta_\pm-1}{\pm\sqrt{(1-\lambda\tau+\beta)^2-4\beta}}.
    \end{equation}
    Note that we have already established that within the considered parameter ranges $|\eta_\pm|\leq \rho(Q)<1$ so that the enumerator in~\eqref{eq:heavy_ball_derivatives} is always negative.
    Thus, in the case $\lambda\tau+\beta>0$ we have $r(\lambda,\beta)=\eta_+>0$ and $\eta_+'<0$ and if in the case $\lambda\tau+\beta<0$ we have $r(\lambda,\beta)=-\eta_->0$ and $\eta_-'>0$. This implies $r(\lambda,\beta)$ is monotonically decreasing in $\beta$.

    Combining both cases we find that $\beta\mapsto r(\lambda,\beta,\tau)$ is decreasing for $\beta<\hat\beta(\lambda)$ and increasing for $\beta>\hat\beta(\lambda)$ and, therefore, minimal for $\beta= \hat\beta(\lambda)$. We will now conclude that
    \begin{equation}
        \beta^*(\tau) = \max_{i}\hat\beta(\lambda_i).
    \end{equation}
    Indeed, note that for all $\lambda_i$, we have $r(\lambda_i,\beta)=\sqrt{\beta}$ and, thus, also $\rho(\tau, \beta) = \sqrt{\beta}$ for $\beta\in [\beta^*,1)$. On the other hand-side for $\beta<\beta^*$ we may pick $i_0$ such that $\beta^* = \hat\beta(\lambda_{i_0})$. Then, $r(\lambda_{i_0},\beta)>r(\lambda_{i_0},\beta^*) = \sqrt{\beta^*}$ as $r(\lambda_{i_0},\beta)$ is decreasing for $\beta<\hat\beta(\lambda_{i_0})$ which implies that $q(\tau,\beta)>\sqrt{\beta^*}$ also for $\beta<\beta^*$ and the result follows. To conclude the proof we note that since $\lambda\mapsto(1-\sqrt{\lambda\tau})^2$ is convex,
    \begin{equation}
        \max_i \hat\beta(\lambda_i) = \max_{\lambda\in\{m,L\}}  (1-\sqrt{\lambda\tau})^2.
    \end{equation}
\end{proof}

\subsubsection{Convergence rate}
\begin{lemma}\label{lemma_Q_power}
    Define the spectral gap
    \begin{equation}
        \varpi \coloneqq \min_{\lambda} \lambda-m,\quad \text{s.t. $\lambda$ is an eigenvalue of $\Lambda$ and $\lambda>m$}
    \end{equation}
    \ie, the difference between the smallest and the second-smallest eigenvalue. Assume the step size satisfies $\tau\leq \frac{1}{L}$. Then with the optimal momentum $\beta=\beta^*(\tau)$, the matrix $Q = Q(\tau,\beta^*(\tau))$ satisfies
    \begin{equation}
        \|Q^k\|
        \lesssim \frac{1}{\varpi (1-\frac{1}{\sqrt{\kappa}})}\rho(Q)^{k-1} \sqrt{2+k^2}
    \end{equation}
    if $\kappa>1$ and
    \begin{equation}
        \|Q^k\|
        \lesssim (1-\sqrt{m\tau})^{k-1} \sqrt{2+k^2}
    \end{equation}
    if $\kappa=1$.
\end{lemma}
\begin{proof}
    We begin by deriving an explicit form of the matrix $P$ which transforms $Q$ into its Jordan normal form, that is
    \begin{equation}
        P^{-1} Q P = \begin{pmatrix}
        J_1 & & \\
        & \ddots & \\
        & & J_k
        \end{pmatrix},
        \qquad
        J_i =
        \begin{pmatrix}
            \mu_i & 1 & & \\
            & \mu_i & \ddots & \\
            & & \ddots & 1 \\
            & & & \mu_i.
        \end{pmatrix}
    \end{equation}
    with $\mu_i$ the eigenvalues of $Q$ which are all of the form $\eta_\pm$ from~\eqref{eq:eta}. Due to the specific block structure of $Q$, we can find $P$ by deriving the Jordan normal form of all $2\times 2$ blocks of the form
    \begin{equation}
        M_i = \begin{pmatrix}
            1-\lambda_i\tau + \beta & -\beta\\
            1& 0
        \end{pmatrix}.
    \end{equation}
    We make a case distinction. If $\eta_+ \neq \eta_+$, $M_i$ is diagonalizable with the matrix
    \begin{equation}
        P_i = \begin{pmatrix}
            \eta_+ & \eta_-\\
            1&1
        \end{pmatrix}.
    \end{equation}
    Conversely, if $\eta_+ = \eta_- = \sqrt{\beta}$, $M_i$ admits a Jordan block and the transformation matrix $P_i$ reads as
    \begin{equation}
        P_i = \begin{pmatrix}
            \sqrt{\beta}& 1\\
            1& 0
        \end{pmatrix}.
    \end{equation}
    The full transformation $P$ can be built from the $P_i$ using Kronecker products and we have $\|Q^k\|\leq \|P\|\|P^{-1}\| \|J^k\|$. 
    Moreover, we have for the spectral norms $\|J^k\|\leq \max_i\|J_i^k\|$ and similarly for $P_i$ and $P_i^{-1}$.
    For the norms of $P_i$ and $P^{-1}_i$ in the case $\eta_+\neq \eta_-$ we find for the Frobenius norms $\|P_i\|^2\leq \|P_i\|_F^2 = 2+|\eta_+|^2 + |\eta_-|^2 \leq 4$ since $|\eta_\pm|\leq 1$ and 
    \begin{equation}
        \begin{aligned}
            \|P^{-1}_i\|^2\leq\|P^{-1}_i\|_F^2
            \leq \frac{4}{\det(P_i)^2}
            = \frac{4}{|\eta_+-\eta_-|^2}.
        \end{aligned}
    \end{equation}
    Moreover, 
    \begin{equation}
        \begin{aligned}
            |\eta_+-\eta_-|^2 
            =& (\lambda-\friction^2)(4(1-\friction\sqrt{\tau})-\tau(\lambda-\friction^2))\\
            =& (\lambda-\friction^2) (4-\sqrt{\tau} (4\friction + \sqrt{\tau}(\lambda-\friction^2)))\\
            \geq& (\lambda-m) \frac{2\sqrt{L}(\sqrt{L}-\sqrt{m})+(\sqrt{L}-\sqrt{m})^2}{L}\\
            \geq& (\lambda-m) \frac{(\sqrt{L}-\sqrt{m})(3\sqrt{L}-\sqrt{m})}{L}\\
            \geq& 2(\lambda-m) (1-\frac{1}{\sqrt{\kappa}}).
        \end{aligned}
    \end{equation}
    In the case $\eta_+=\eta_-$ we simply have $\|P\|^2 \leq 3$ and $\|P^{-1}\|^2 \leq 3$.
    It is well known that for one Jordan block we find
    \begin{equation}
        J_i^k = \begin{pmatrix}
            \mu_i^k & k\mu_i^{k-1} & \frac{k(k-1)\mu_i^{k-2}}{2}& \dots& \binom{k}{j}\mu_i^{k-j}&\dots\\
            & \mu_i^k & k\mu_i^{k-1} &\ddots & &\\
            & & \ddots & \ddots & & \\
            & & & & & \\
            & & & & & \\
            & & & & & \mu_i^k.
        \end{pmatrix}
    \end{equation}
    Since in our case each Jordan block is of size at most two we find
    \begin{equation}
        \begin{aligned}
            \|J_i^k\|^2 \leq \|J_i^k\|^2_F 
            \leq 2|\mu_i|^{2k} + k^2 |\mu_i|^{2(k-1)}
            \leq \rho(Q)^{2(k-1)} (2\rho(Q)^2 + k^2)
            \leq \rho(Q)^{2(k-1)} (2+k^2).
        \end{aligned}
    \end{equation}
    In total we therefore find
    \begin{equation}
        \|Q^k\|\leq \|P\|\|P^{-1}\| \rho(Q)^{k-1} \sqrt{2+k^2}
        \lesssim \frac{1}{\varpi (1-\frac{1}{\sqrt{\kappa}})}\rho(Q)^{k-1} \sqrt{2+k^2}
    \end{equation}
    which implies the desired result upon noting that for $\tau<\frac{1}{L}$, $\rho(Q) = (1-\sqrt{m\tau})$.
\end{proof}

\begin{corollary}\label{convergence_corollary}
    For every $\delta\in (0,1-\rho(Q))$ there exists a $K>0$ with
    \begin{equation}
        K = \tilde\Oc(\delta^{-1})
    \end{equation}
    where $\tilde\Oc$ hides logarithmic factors, such that for any $k\geq K$
    \begin{equation}
        \|Q^k\|
        \lesssim  ((1-\sqrt{m\tau})+\delta)^k.
    \end{equation}
\end{corollary}
\begin{proof}
    We consider only the case $\kappa>1$. The converse case is handled similarly. By~\cref{lemma_Q_power} we have
    \begin{equation}
        \|Q^k\|^{\frac{1}{k-1}}
        \lesssim (1-\sqrt{m\tau}) \left(\frac{\sqrt{2+k^2}}{\varpi(1-\frac{1}{\sqrt{\kappa}})}\right)^{\frac{1}{k-1}}
        \lesssim (1-\sqrt{m\tau}) \left(\frac{k}{\varpi (1-\frac{1}{\sqrt{\kappa}})}\right)^{\frac{1}{k-1}}.
    \end{equation}
    By~\cref{lemma:poly_bound_by_exp} it holds
    \begin{equation}
        k\leq \frac{4}{\log(2)\delta}(1+\delta/2)^{k-1}
    \end{equation}
    for all $k\geq 0$. Thus, we have
    \begin{equation}
        (1-\sqrt{m\tau}) \left(\frac{k}{\varpi(1-\frac{1}{\sqrt{\kappa}})}\right)^{1/(k-1)}
        \lesssim(1+\delta/2)(1-\sqrt{m\tau}) \left(\frac{1}{\delta \varpi (1-\frac{1}{\sqrt{\kappa}})}\right)^{1/(k-1)}.
    \end{equation}
    For
    \begin{equation}\label{eq_k:ploylog}
        k-1\gtrsim\frac{\log\left(\frac{1}{\delta \varpi (1-\frac{1}{\sqrt{\kappa}})}\right)}{\log\left(1+\frac{\delta}{2}\right)}.
    \end{equation}
    it holds 
    \begin{equation}
        (1+\delta/2)\left(\frac{1}{\delta \varpi (1-\frac{1}{\sqrt{\kappa}})}\right)^{1/(k-1)}
        \leq (1+\delta)
        \leq \bigg(1+\frac{\delta}{1-\sqrt{m\tau}}\bigg).
    \end{equation}
    Thus, for such $k$
    \begin{equation}
        \|Q^k\|
        \leq (\|Q^k\|^{1/(k-1)})^{k-1}
        \leq
        \left((1-\sqrt{m\tau})(1+\frac{\delta}{1-\sqrt{m\tau}})\right)^{k-1}
        \leq(1-\sqrt{m\tau}+\delta)^{k-1}
    \end{equation}
    as desired. It remains te estimate the scaling of $K$. 
    Since for $\delta\in (0,1)$, $\log(1+\delta/2)\gtrsim \delta$ we find
    \begin{equation}
        \frac{\log\left(\frac{1}{\delta \varpi (1-\frac{1}{\sqrt{\kappa}})}\right)}{\log\left(1+\frac{\delta}{2}\right)}
        \leq \delta^{-1}\log\left(\frac{1}{\delta \varpi (1-\frac{1}{\sqrt{\kappa}})}\right).
    \end{equation}
    Moreover, note that for some eigenvalue $\lambda$ of $\Lambda$
    \begin{equation}
        \begin{aligned}
            \varpi (1-\frac{1}{\sqrt{\kappa}})
            = (\lambda-m)\frac{\sqrt{L}-\sqrt{m}}{\sqrt{L}}
            &= (\lambda-m)\frac{L-m}{\sqrt{L}(\sqrt{L}+\sqrt{m})}\geq \frac{\varpi^2}{\sqrt{L}(\sqrt{L}+\sqrt{m})} \geq \frac{\varpi^2}{2L}
        \end{aligned}
    \end{equation}
    \begin{equation}
        \frac{\log\left(\frac{1}{\delta \varpi (1-\frac{1}{\sqrt{\kappa}})}\right)}{\log\left(1+\frac{\delta}{2}\right)}
        \leq \delta^{-1}\log\left(\frac{2L}{\delta \varpi^2}\right).
    \end{equation}
\end{proof}

\begin{lemma}\label{lemma:poly_bound_by_exp}
    For any $\delta\in (0,1)$ it holds
    \begin{equation}
        k\leq \frac{2}{\log(2)\delta}(1+\delta)^{k-1},\quad k\geq 1
    \end{equation}
    and for any $q\in (0,1)$, $\delta\in (0,1-q)$ it holds
    \begin{equation}
        k q^{k-1} \leq \big(\frac{1}{q} + \frac{1}{\delta}\big)(q+\delta)^k,\quad k\geq 1.
    \end{equation}
\end{lemma}
\begin{proof}
    Our goal is to find a constant $c$ such that for all $k$, $k\leq c(1+\delta)^{k-1}$. Indeed, we have
    \begin{equation}
        c\leq \max_{x\geq 1}  \frac{x}{(1+\delta)^{x-1}}
    \end{equation}
    or, equivalently
    \begin{equation}
        \log(c) \leq \max_{x\geq 1}  \log(x) - (x-1)\log(1+\delta).
    \end{equation}
    The right-hand side is increasing for $x \leq \frac{1}{\log(1+\delta)}$ and decreasing for $x \geq \frac{1}{\log(1+\delta)}$ and, thus, maximal for $x = \frac{1}{\log(1+\delta)}$ since $\delta<1$ and $\log(2)<1$. Therefore, 
    \begin{equation}
        c(\delta) \leq \frac{1}{\log(1+\delta)(1+\delta)^{\frac{1}{\log(1+\delta)}-1}}
    \end{equation}
    and since $(1+\delta)^{\frac{1}{\log(1+\delta)}}\geq 1$
    \begin{equation}
        c(\delta) \leq \frac{1+\delta}{\log(1+\delta)}.
    \end{equation}
    Lastly, by concavity of the logarithm and the assumption $\delta\in (0,1)$ we have $\log(1+\delta)>\log(2)\delta$ yielding
    \begin{equation}
        c(\delta)\leq \frac{2}{\log(2)\delta}.
    \end{equation}
    For the second assertion we proceed similarly. We need to compute
    \begin{equation}\label{eq:polybound}
        \max_{x} \frac{xq^{x-1}}{(q+\delta)^x}
    \end{equation}
    which, after applying a logarithm reads as
    \begin{equation}
        \max \log(x) + (x-1)\log(q) - x\log(q+\delta).
    \end{equation}
    The maximum is attained for 
    \begin{equation}
        x = \log\big(\frac{q+\delta}{q}\big)^{-1}
    \end{equation}
    if this quantity is greater than $1$ and for $x=1$ else. In the latter case we obtain
    \begin{equation}
        \max_{x} \frac{xq^{x-1}}{(q+\delta)^x} = \frac{1}{q+\delta}\leq \frac{1}{q} + \frac{1}{\delta}.
    \end{equation}
    In the former case, inserting into~\cref{eq:polybound} yields
    \begin{equation}\label{eq:polybound2}
        \begin{aligned}
            \max_{x} \frac{xq^{x-1}}{(q+\delta)^x} 
            =&\frac{\log\big(\frac{q+\delta}{q}\big)^{-1}q^{\log\big(\frac{q+\delta}{q}\big)^{-1}-1}}{(q+\delta)^{\log\big(\frac{q+\delta}{q}\big)^{-1}}}\\
            =&\log\big(1+\frac{\delta}{q}\big)^{-1}\frac{1}{q}\bigg(\frac{q}{q+\delta} \bigg)^{\log\big(\frac{q+\delta}{q}\big)^{-1}}\\
            \leq&\log\big(1+\frac{\delta}{q}\big)^{-1}\frac{1}{q}.
        \end{aligned}
    \end{equation}
    By concavity of the logarithm we have
    \begin{equation}
        \log(1+\delta/q) = \log(q+\delta) - \log(q) \geq \frac{\delta}{q+\delta}
    \end{equation}
    Inserting into~\cref{eq:polybound2} finally yields
    \begin{equation}
        \begin{aligned}
            \max_{x} \frac{xq^{x-1}}{(q+\delta)^x} 
            \leq&\frac{q+\delta}{q\delta} = \frac{1}{q} + \frac{1}{\delta}.
        \end{aligned}
    \end{equation}
\end{proof}

\subsection{The General Setting}
\subsubsection{Discrete time}

\begin{lemma}\label{lemma:ila2-recurrence}
    Assume that $\tau,\beta>0$ satisfy $1-L\tau-\beta^2>0$.
    Define the function
    \begin{equation}
        \omega \colon \R^d \times \R^d \to \R, \quad \omega (x, v) = U(x) + c_1\|v\|^2 + c_2\|x\|^2 + c_3\inner{x}{v}
    \end{equation}
    with coefficients
    \begin{equation}\label{eq:discrete_lyapunov_ci}
        c_1 = \frac{1-L\tau+\beta^2}{4}, \quad
        c_2 = \frac{(1-\beta)(1-L\tau-\beta^2)}{4(1+\beta)\tau}, \quad
        c_3 = \beta\frac{1-L\tau-\beta^2}{2(1+\beta)\sqrt{\tau}}.
    \end{equation}
    Then, it holds true that
    \begin{equation}\label{eq:discrete_lyapunov}
        \begin{aligned}
            R_{\tau,\beta}\omega(x,v)\leq& \omega(x,v)
            -\frac{1-L\tau-\beta^2}{4} \|v\|^2
            -\beta\frac{1-L\tau-\beta^2}{2(1+\beta)} \inner{x}{\nabla U(x)}
            +d(1-\beta).
        \end{aligned}
    \end{equation}
\end{lemma}
\begin{proof}
    Consider the general function $\omega(x,v) = U(x) + c_1\|v\|^2 + c_2\|x\|^2 + c_3\inner{x}{v}$ with coefficients $c_1, c_2, c_3 \in \R$.
    We will show that the choices for the coefficients naturally appear in order to obtain the desired drift condition.
    
    By $L$-smoothness of $U$ we have for any $p,q \in \R^d$, $U(q)\leq U(p) + \langle\nabla U(p),q-p\rangle + \frac{L}{2}\|p-q\|^2$. Therefore, we find after performing one \cref{M:eq:ila-2} step starting from $X^0=x$, $V^0=v$ that
    \begin{equation}
        \begin{aligned}
            R_{\tau,\beta}\omega(x,v)\leq& \E\bigg[U(x) + \inner{\nabla U(x)}{\sqrt{\tau} V^1} + \frac{L}{2}\|\sqrt\tau V^1\|^2\\
            &+c_1\|V^1\|^2
            +c_2\|x+\sqrt\tau V^1\|^2
            +c_3\inner{x+\sqrt\tau V^1}{V^1}\bigg]\\
            =&U(x) + c_1\|v\|^2 + c_2\|x\|^2+c_3\inner{x}{v}\\
            &+\inner{x}{v}\left\{ 2\sqrt\tau c_2 \beta + c_3 \beta - c_3  \right\}\\
            &+\inner{v}{\nabla U(x)}\left\{\sqrt{\tau}\beta -2\sqrt\tau\beta\left(\frac{L\tau}{2}+c_1+c_2\tau+c_3\sqrt\tau\right)\right\}\\
            &+\|v\|^2\left\{\beta^2\left(\frac{L\tau}{2}+c_1+c_2\tau+c_3\sqrt\tau\right)-c_1\right\}\\
            &+\inner{x}{\nabla U(x)}\left\{-\sqrt\tau(2c_2\sqrt\tau+c_3)\right\}\\
            &+\|\nabla U(x)\|^2\bigg\{-\tau +\tau \left(\frac{L\tau}{2}+c_1+c_2\tau+c_3\sqrt\tau\right)\bigg\}\\
            &+\left(\frac{L\tau}{2}+c_1+c_2\tau+c_3\sqrt\tau\right)2d(1-\beta).
        \end{aligned}
    \end{equation}
    We ensure that the scalar products $\inner{x}{v}$ and $\inner{v}{\nabla U(x)}$ vanish by setting
    \begin{equation}
            \frac{L\tau}{2}+c_1+c_2\tau+c_3\sqrt\tau = \frac{1}{2}, \quad
             2\sqrt\tau c_2 \beta + c_3 \beta - c_3  = 0
    \end{equation}
    implying that we can express $c_1$ and $c_2$ in terms of $c_3$, to obtain
    \begin{equation}\label{eq:c2-inequalities}
        c_1 = \frac{1-L\tau -c_3\sqrt{\tau}\frac{1+\beta}{\beta}}{2}, \quad c_2 = \frac{(1-\beta)c_3}{2\beta\sqrt\tau}.
    \end{equation}
    Moreover, these conditions lead to 
    \begin{equation}\label{eq:discrete_lyapunov_proof}
        \begin{aligned}
            R_{\tau,\beta}\omega(x,v)\leq&\omega (x, v)\\
            &\hspace{-2cm}+\left\{\frac{\beta^2}{2}-c_1\right\}\|v\|^2+\left\{-\sqrt\tau(2c_2\sqrt\tau+c_3)\right\}\inner{x}{\nabla U(x)} -\frac{\tau}{2}\|\nabla U(x)\|^2 +d(1-\beta),
        \end{aligned}
    \end{equation}
    from which we can see that we require
    \begin{equation}\label{eq:c3-inequalities}
            \frac{\beta^2}{2}-c_1 <0, \quad-2c_2\sqrt\tau - c_3 <0.
    \end{equation}
    Inserting~\cref{eq:c2-inequalities} and noting that the constants should be stable as $\tau\rightarrow 0$ we arrive at
    \begin{equation}
        0 < c_3 < \frac{(1 - L\tau - \beta^2)\beta}{\sqrt{\tau} (1 + \beta)},
    \end{equation}
    which is achieved by our choice
    \begin{equation}
        c_3 = \frac{(1 - L\tau - \beta^2)\beta}{2\sqrt{\tau} (1 + \beta)}
    \end{equation}
    resulting in~\eqref{eq:discrete_lyapunov_ci}.
    Inserting the $c_i$ into~\eqref{eq:discrete_lyapunov_proof} and omitting the remaining non-positive $\Oc(\tau)$ terms, which will be dominated by $\Oc(\sqrt{\tau})$ for small $\tau$, we obtain the desired result.
    \end{proof}

\begin{lemma}\label{lemma_lyapunov_equivalent}
        Assume that either $\beta$ is fixed, or satisfies~\Cref{M:ass:beta}.
        With the choices for $c_1, c_2, c_3$ from \cref{M:lemma:ila2-recurrence} there exist $\xi_1, \xi_2,\xi_3, \xi_4>0$ such that for any $\tau>0$ with $1-L\tau-\beta^2>0$ it holds true that
        \begin{equation}\label{eq:lower bound of lyapunov}
            U(x) +\zeta_1 \|x\|^2 + \zeta_2\|v\|^2\leq \omega(x,v)\leq U(x) +\zeta_3 \|x\|^2 + \zeta_4\|v\|^2
        \end{equation}
\end{lemma}
\begin{proof}
        By Young's inequality
        we have for any $\delta>0$
        \begin{equation}
            \omega(x, v) 
            \geq U(x) + \left(c_1-c_3\frac{\delta}{2}\right)\|v\|^2 + \left(c_2-c_3\frac{1}{2\delta}\right)\|x\|^2.
        \end{equation}
        Consequently, it is enough to show that there exists $\delta>0$ satisfying $c_1-c_3\frac{\delta}{2} > 0$ and $c_2-c_3\frac{1}{2\delta}>0$ and that these inequalities remain strict as $\sqrt{\tau}\rightarrow 0$ under~\Cref{M:ass:beta}.
        Denoting for simplicity $a=1-L\tau -\beta^2$ and inserting choices for $c_1, c_2, c_3$ from \cref{M:lemma:ila2-recurrence} we have to ensure that
        \begin{equation}
                \frac{a+2\beta^2}{4}-\frac{\delta}{2}\frac{a\beta}{2(1+\beta)\sqrt\tau}>0 \quad \text{and} \quad
                \frac{(1-\beta)a}{4(1+\beta)\tau}-\frac{1}{2\delta}\frac{a\beta}{2(1+\beta)\sqrt{\tau}}>0.
        \end{equation}
        Rearranging for $\delta$ leads to the bounds
        \begin{equation}
                \frac{\beta\sqrt{\tau}}{1-\beta} < \delta < \frac{(a+2\beta^2)(1+\beta)\sqrt{\tau}}{a\beta}.
        \end{equation}
        In order to conclude that such a $\delta$ exists, we need to show that, indeed,
        \begin{equation}\label{eq:equivalent lyapunov}
            \frac{\beta\sqrt{\tau}}{1-\beta} < \frac{(a+2\beta^2)(1+\beta)\sqrt{\tau}}{a\beta}.
        \end{equation}
        Inserting and simplifying shows that this condition is equivalent to $0<(1-L\tau) - \beta^2(1-2L\tau)$.
        It is obvious that the above condition is satisfied for $\frac{1}{2L}\leq \tau <\frac{1}{L}$. On the other hand-side, for $\tau<\frac{1}{2L}$ we can divide by $(1-2L\tau)$ and obtain
        \begin{equation}
            \beta^2<\frac{1-L\tau}{1-2L\tau} = 1 + \frac{L\tau}{1-2L\tau},
        \end{equation}
        which is again satisfied for any $\beta\in(0,1)$. Moreover, one may check that 
        \begin{equation}
            \delta = \frac{(\beta + \frac{1}{2\beta})\sqrt{\tau}}{1-\beta}
        \end{equation}
        satisfies the constraints for all $\tau\in [0,\frac{1}{L})$ and yields
        \begin{equation}
            \begin{aligned}
                c_1-c_3\frac{\delta}{2}\rightarrow \frac{2-\beta}{8}>\frac{1}{8},\quad c_2-c_3\frac{1}{2\delta}\rightarrow +\infty
            \end{aligned}
        \end{equation}
        as $\tau\rightarrow0$ which implies the existence of $\xi_1,\xi_2>0$ independent of $\tau>0$ such that~\eqref{eq:lower bound of lyapunov} holds.
        
        Conversely, under~\Cref{M:ass:beta} with the same choice of $\delta$ it follows
        \begin{equation}
            \begin{aligned}
                c_1-c_3\frac{\delta}{2}\rightarrow \frac{1}{8},\quad c_2-c_3\frac{1}{2\delta}\rightarrow \frac{\friction^2}{3}
            \end{aligned}
        \end{equation}
        so that also in this case $\xi_1, \xi_2$ can be chosen uniformly over $\tau>0$. The upper bound of $\omega$ follows similarly.
\end{proof}

\begin{proposition}\label{prop:lyapunov_discrete}
    Let~\Cref{M:ass:assumptions} hold and assume that $\tau,\beta>0$ satisfy $1-L\tau-\beta^2>0$.
    Consider $\omega \colon \R^d \times \R^d \to \R$ from \cref{M:lemma:ila2-recurrence}.
    Then there exist $\eta,b>0$, $\lambda\in (0,1)$ explicit in the proof and all independent of $\tau$ and $\beta$ such that for $\Vc(x,v)\coloneqq\expo{\eta\sqrt{1+\omega(x,v)}}$ it holds that
    \begin{equation}\label{eq:lyapunov-conditioon}
        R_{\tau,\beta}\Vc(x,v)\leq \lambda^{1-\beta}\Vc(x,v) + (1-\beta)b.
    \end{equation}
\end{proposition}
\begin{proof}
    As in the proof of~\cref{M:prop:lyapunov_stronger} we obtain
    \begin{equation}\label{eq:lyapunov_discrete_final1}
        \begin{aligned}
            R_{\tau,\beta}\Vc(x,v)
            &\leq \expo{\frac{\eta^2(1-\beta)M^2}{2}} \expo{\eta\bigl(1+\omega(x,v)
            -\frac{1-L\tau-\beta^2}{4} \|v\|^2
            -\beta\frac{1-L\tau-\beta^2}{2(1+\beta)} \inner{x}{\nabla U(x)}
            +d(1-\beta)\bigr)^{1/2}}.
        \end{aligned}
    \end{equation}
    ans again by the same arguments as in~\cref{M:prop:lyapunov_stronger},
    \begin{equation}\label{eq:lyapunov_discrete_final2}
        \begin{aligned}
            \biggl(1+\omega(x,v)&
            -\frac{1-L\tau-\beta^2}{4} \|v\|^2
            -\beta\frac{1-L\tau-\beta^2}{2(1+\beta)} \inner{x}{\nabla U(x)}
            +d(1-\beta)\biggr)^{1/2}\\
            \leq& (1+\omega(x,v))^{1/2} - (1-\beta)\frac{\|v\|^2\frac{1-L\tau-\beta^2}{4(1-\beta)}+\inner{x}{\nabla U(x)}\beta\frac{1-L\tau-\beta^2}{2(1-\beta^2)}-d}{2(1+\omega(x,v))^{1/2}}
        \end{aligned}
    \end{equation}
    By \Cref{M:ass:assumptions}, there exist $R_x,R_v>0$ which can\footnote{contrary to~\cref{M:prop:lyapunov_stronger}} be chosen uniformly with respect to $\tau$ such that $R_x\sim d^{1/\grc}$, $R_v\sim d^{1/2}$, and $\rho>0$ such that for $\|x\|\geq R_x$ or $\|v\|\geq R_v$ it holds 
    \begin{equation}
        \frac{\|v\|^2\frac{1-L\tau-\beta^2}{4(1-\beta)}+\inner{x}{\nabla U(x)}\beta\frac{1-L\tau-\beta^2}{2(1-\beta^2)}-d}{2(1+\omega(x,v))^{1/2}}\geq \rho
    \end{equation}
    Thus, with $\eta=\rho/M^2$ and $\lambda \coloneqq \expo{-\frac{\rho^2}{2M^2}}\in (0,1)$ we find
    \begin{equation}
        \begin{aligned}
            R_{\tau,\beta}\Vc(x,v)
            \leq& \lambda^{1-\beta}\Vc(x,v).
        \end{aligned}
    \end{equation}
    On the other hand, for $\|x\|\leq R_x$ and $\|v\|\leq R_v$ we simply estimate from~\cref{eq:lyapunov_discrete_final1}, \eqref{eq:lyapunov_discrete_final2}
    \begin{equation}\label{eq:lyapunov_discrete_final3}
        \begin{aligned}
            R_{\tau,\beta}\Vc(x,v)\leq& \lambda^{1-\beta}\Vc(x,v) \\
            &+ \lambda^{1-\beta}\Vc(x,v)\Biggl( \expo{(1-\beta)\biggl(\frac{\rho^2}{2M^2}-\frac{\|v\|^2\frac{1-L\tau-\beta^2}{4(1-\beta)}+\inner{x}{\nabla U(x)}\beta\frac{1-L\tau-\beta^2}{2(1-\beta^2)}-d}{2(1+\omega(x,v))^{1/2}}\biggr)} - 1\Biggr).
        \end{aligned}
    \end{equation}
    By convexity of the exponential, it holds for any $t>0$, $\expo{t}-1\leq \expo{t}t$ and for $t\leq 0$ we have $\expo{t}-1\leq 0\leq \expo{t}|t|$. Thus, denoting  
    \begin{equation}
        C(x,v,\tau) = 
        \frac{\rho^2}{2M^2}-\frac{\|v\|^2\frac{1-L\tau-\beta^2}{4(1-\beta)}+\inner{x}{\nabla U(x)}\beta\frac{1-L\tau-\beta^2}{2(1-\beta^2)}-d}{2(1+\omega(x,v))^{1/2}}<\infty
    \end{equation}
    we obtain from~\eqref{eq:lyapunov_discrete_final3}
    \begin{equation}
        \begin{aligned}
            R_{\tau,\beta}\Vc(x,v)
            \leq& \lambda^{1-\beta}\Vc(x,v)+ \lambda^{1-\beta}\Vc(x,v)\expo{(1-\beta)C(x,v,\tau)}(1-\beta)|C(x,v,\tau)|\\
            \leq& \lambda^{1-\beta}\Vc(x,v)+ (1-\beta)\sup_{\substack{\tau\in (0,1/L)\\ \|x\|\leq R_x,\,\|v\|\leq R_v}}\lambda^{1-\beta}\Vc(x,v)\expo{(1-\beta)C(x,v,\tau)}|C(x,v,\tau)|\\
        \end{aligned}
    \end{equation}
    concluding the proof with $b$ the supremum above whose 
    finiteness follows from continuity of the involved functions in $(x,v)$.
\end{proof}

\subsubsection{Continuous time}
\begin{proposition}\label{prop:lyapunov_cont}
    Let $\Vc$ as in~\cref{M:prop:lyapunov_discrete}, albeit with the parameters
    \begin{equation}
        c_1 = \frac{1}{2},\quad c_2 = \frac{\friction^2}{4},\quad c_3 = \frac{\friction}{2}.
    \end{equation}
    Then there exist $\lambda,b>0$ which are explicit in the proof such that
    \begin{equation}\label{eq:lyapunov_drift_cont}
        \Lc\Vc(x,v)\leq -\lambda \Vc(x,v) + b.
    \end{equation}
\end{proposition}
\begin{proof}
    We compute
    \begin{equation}
        \begin{aligned}
            \frac{\Lc \Vc(x,v)}{\Vc(x,v)}
            =& \frac{\eta}{2\sqrt{1+\omega(x,v)}}\bigg\{
            \inner{x}{v}\big(2c_2-\friction c_3\big)\\
            &+\inner{\nabla U(x)}{v}\big(1-2c_1\big)\\
            &+\|v\|^2\big(c_3-2\friction c_1\big)\\
            &+\inner{\nabla U(x)}{x}\big(-c_3\big)\\
            &+\friction\eta\frac{\|2c_1v+c_3x\|^2}{4\sqrt{1+\omega(x,v)}}
            + 2\friction c_1 d
            -\friction\frac{\|2c_1v+c_3x\|^2}{4(1+\omega(x, v))}
            \bigg\}.
        \end{aligned}
    \end{equation}
    With the stated choices for the coefficients, it follows that
    \begin{equation}\label{eq:lyap_cont1}
        \begin{aligned}
            \frac{\Lc \Vc(x,v)}{\Vc(x,v)}
            \leq & \frac{\eta}{2\sqrt{1+\omega(x,v)}}\bigg\{
            -\frac{\friction}{2}\|v\|^2
            -\frac{\friction}{2}\inner{\nabla U(x)}{x}
            +\friction\eta\frac{\|2c_1v+c_3x\|^2}{4\sqrt{1+\omega(x,v)}}
            + 2\friction c_1 d
            \bigg\}.
        \end{aligned}
    \end{equation}
    By~\Cref{M:ass:assumptions} and since
    \begin{equation}
        \frac{\|2c_1v+c_3x\|^2}{4\sqrt{1+\omega(x,v)}}\sim \|(x,v)\|
    \end{equation}
    we can find $R_x, R_v,\rho,\eta>0$ with $R_v \sim d^{1/2}$ and $R_x \sim d^{1/\grc}$ such that for $\|x\|\geq R_x$ or $\|v\|\geq R_v$ we have
    \begin{equation}
        \frac{\friction}{2}\|v\|^2
        +\frac{\friction}{2}\inner{\nabla U(x)}{x}-\friction\eta\frac{\|2c_1v+c_3x\|^2}{4\sqrt{1+\omega(x,v)}} - 2\friction c_1 d\geq \rho (\|x\| + \|v\|^2).
    \end{equation}
    Thus, we obtain for such $(x,v)$
    \begin{equation}
        \begin{aligned}
            \frac{\Lc \Vc(x,v)}{\Vc(x,v)}
            \leq & -\frac{\eta \rho (\|x\| + \|v\|^2)}{2\sqrt{1+\omega(x,v)}}\leq -\lambda
        \end{aligned}
    \end{equation}
    for some $\lambda>0$ since $\sqrt{1+\omega(x,v)}\sim \|(x,v)\|$.
    On the other hand-side, for $\|x\|\leq R_x$ and $\|v\|\leq R_v$ we have from~\eqref{eq:lyap_cont1}
    \begin{equation}
        \begin{aligned}
            \Lc \Vc(x,v)
            \leq &-\lambda \Vc(x,v) \\
            &+ \Vc(x,v) \bigg(\frac{\eta}{2\sqrt{1+\omega(x,v)}}\bigg\{
            \frac{\epsilon}{2}\|\nabla U(x)\|\|x\|
            +\friction\eta\frac{\|2c_1v+c_3x\|^2}{4\sqrt{1+\omega(x,v)}}
            + 2\friction c_1 d
            \bigg\} - \lambda\bigg)
        \end{aligned}
    \end{equation}
    so that the result follows with 
    \begin{equation}
        b \coloneqq \sup_{\|(x,v)\|\leq R}\Vc(x,v) \bigg(\frac{\eta}{2\sqrt{1+\omega(x,v)}}\bigg\{
            \frac{\epsilon}{2}\|\nabla U(x)\|\|x\|
            +\friction\eta\frac{\|2c_1v+c_3x\|^2}{2\sqrt{1+\omega(x,v)}}
            + 2\friction c_1 d
            \bigg\} - \lambda\bigg).
    \end{equation}
\end{proof}

\begin{proposition}[Discretization error]\label{thm:disc_SM}
    Let~\Cref{M:ass:assumptions,M:ass:beta} hold true and denote $\step\coloneqq\sqrt{\tau}$. Then there exists $\rho, c >0$ such that for any initial distribution $\mu\in\Pc_2(\R^{2d})$ it holds true
    \begin{equation}
        \Wc_2^2(\mu R_{\tau,\beta}^k,\mu P_{k\step})\leq k \disca \step^2\expo{c k\step},\quad k\geq 0.
    \end{equation}
\end{proposition}
\begin{proof}
    We first analyze the one-step discretization error.
    To do so, we require a representation of~\eqref{M:eq:ila-2} as an adapted stochastic process in continuous time, namely, $\b{Z}_t = (\b{X}_t,\b{V}_t)$ defined piecewise via
    \begin{equation}\label{eq:disc_sde_interp}
        \begin{cases}
            \dd \b{X}_t &= \left[ \b{V}_t + (t-t_k)(-2\friction \b{V}_{t_k} - \nabla U(\b{X}_{t_k})) \right]\dd t + (t-t_k)\sqrt{4\friction}\dd W_t,\\
            \dd \b{V}_t &= \left[-2\friction \b{V}_{t_k} - \nabla U(\b{X}_{t_k})\right]\dd t + \sqrt{4\friction}\dd W_t,
        \end{cases}
    \end{equation}
    for $t\in(t_k,t_{k+1}]$. The second line is easily seen to be an interpolation of~\eqref{M:eq:ila}.
    For the first line note that the process $\b{X}_t = \b{X}_{t_k} + (t-t_k)\b{V}_t$ is adapted and an interpolation of~\eqref{M:eq:ila}.
    Applying \Ito{}'s formula to this process leads to \cref{eq:disc_sde_interp}. We consider in the following for simplicity the step from $k=0$ to $k=1$.
    We begin with the velocity component. 
    Using synchronous coupling by \Ito{}'s lemma we have for any $h\in(0,\step)$
    \begin{equation}\label{eq:discretization_convergence1}
        \begin{aligned}
            \|V_{h} - \b{V}_{h}\|^2 - \|V_{0} - \b{V}_{0}\|^2
            =&2\int_{0}^{h} \begin{bmatrix}
                -\friction \b{V}_{0} - \nabla U(\b{X}_{0})\\
                -\friction V_t - \nabla U(X_t)
            \end{bmatrix}\cdot \begin{bmatrix}
                \b{V}_t-V_t\\
                V_t - \b{V}_t
            \end{bmatrix}\dd t\\
            =& 2\int_{0}^{h}\bigl( \friction V_t -\friction \b{V}_{t} + \friction \b{V}_{t}-\friction \b{V}_{0} \\
            &+\nabla U(X_t) - \nabla U(\b{X}_{t}) + \nabla U(\b{X}_{t}) - \nabla U(\b{X}_{0})\bigr) \bigl( \b{V}_t-V_t\bigr)\dd t\\
            \lesssim& L\int_{0}^{h} \friction^2\|\b{V}_{0} -\b{V}_t\|^2 + \|\b{X}_{0}-\b{X}_t\|^2 + \|\b{X}_t-X_t\|^2 + \|\b{V}_t-V_t\|^2\dd t,
        \end{aligned}
    \end{equation}
    where we used Cauchy-Schwarz, the $L$-smoothness of the potential and Young's inequality.
    Note that $\b{V}_t - \b{V}_{0} = (-\friction \b{V}_{0} - \nabla U(\b{X}_{0}))t + \sqrt{2\friction} (W_t-W_{0})$ and, thus, by the properties of (increments of) Brownian motion, we obtain
    \begin{equation}
        \begin{aligned}
            \E\left[\| \b{V}_t - \b{V}_{0} \|^2\right] 
            &\lesssim t^2\left(\friction^2\E\left[ \|\b{V}_{0}\|^2\right] + \E\left[ L^2\|\b{X}_{0}\|^2\right] \right) + \friction t
        \end{aligned}
    \end{equation}
    where we used that $0 \in \argmin_x U(x)$. Similarly, from the definition of the interpolated position, it follows easily that 
    \begin{equation}
        \begin{aligned}
            \Exp{\|\b{X}_t-\b{X}_{0}\|^2} &= t^2 \Exp{\|\b{V}_t\|^2}   \\
            &\lesssim t^2\left(\E\left[\|\b{V}_{0}\|^2\right] + L^2\E\left[\|\b{X}_{0}\|^2\right]\right)+ \friction t^3.
        \end{aligned}
    \end{equation}
    Taking the expectation of \cref{eq:discretization_convergence1} and inserting these estimates, we obtain
   \begin{equation}\label{eq:v-disc}
        \begin{aligned}
            \Exp{\|\b{V}_{h} - V_{h}\|^2} - \Exp{\|\b{V}_{0} - V_{0}\|^2}&\\
            &\hspace{-6cm}\lesssim h^3 \left((1+\friction^4)\E\left[ \|\b{V}_{0}\|^2\right] + (1+\friction^2)\E\left[ \|\b{X}_{0}\|^2\right] \right) + h^2\friction(1+\friction^2)  \\
            &\hspace{-5cm}+\int_{0}^{h} \Exp{\|\b{X}_t-X_t\|^2} 
            + \Exp{\|\b{V}_t-V_t\|^2}\dd t.
        \end{aligned}
    \end{equation}
    Next, let us consider the analogously coupled $X$ process.
    We apply \Ito{}'s lemma to obtain an \gls{sde} for $\|X_t-\b{X}_t\|^2$.
    While in this case, the diffusion coefficient is nonzero, it nonetheless vanishes in expectation (\cite[Section 3.3]{khasminskii2012stochastic}, \cite[Chapter 8, Section 2]{gikhman1965introduction}).
    Thus, we find
    \begin{equation}\label{eq:x-disc}
        \begin{aligned}
            &\Exp{\|\b{X}_{h} - X_{h}\|^2} - \Exp{\|\b{X}_{0} - X_{0}\|^2}\\
            &= 2\Exp{\int_{0}^{h} \begin{bmatrix}
                \b{V}_{t} + t(-\friction \b{V}_{0} - \nabla U(\b{X}_{0}))\\
                V_t
            \end{bmatrix}\cdot \begin{bmatrix}
                \b{X}_t-X_t\\
                X_t - \b{X}_t
            \end{bmatrix}
            +d\friction t^2 \dd t}\\
            &= 2\int_{0}^{h} \Exp{\left(
                \b{V}_{t} -V_t + t(-\friction \b{V}_{0} - \nabla U(\b{X}_{0})\right)
                \left(\b{X}_t-X_t\right)}\dd t
            +d\friction \frac{h^3}{3}\\
            &= \int_{0}^{h} \Exp{\|\b{V}_t -V_t\|^2 + 2\|\b{X}_t -X_t\|^2 + t^2 \|\friction \b{V}_{0} +\nabla U(\b{X}_{0})\|^2 }\dd t\\
            &\qquad+d\friction \frac{h^3}{3}\\
            &\lesssim \left( \friction^2\E\left[\|\b{V}_{0}\|^2\right] + \E\left[\|\b{X}_{0}\|^2\right] + d\friction\right)h^3 + \int_{0}^{h} \Exp{\|\b{V}_t -V_t\|^2} + \Exp{\|\b{X}_t -X_t\|^2} \dd t.
        \end{aligned}
    \end{equation}
    Combining \cref{eq:v-disc} and \cref{eq:x-disc}, we find an upper bound on the one-step discretization error induced by~\eqref{M:eq:ila-2} as a function of the discretization time. It is given by
    \begin{equation}
        \begin{aligned}
            &\Exp{\|Z_{h} - \b{Z}_{h}\|^2} - \Exp{\|Z_{0} - \b{Z}_{0}\|^2}\\
            &\lesssim \left( (1+\friction^2+\friction^4)\E\left[\|\b{V}_{0}\|^2\right] + (1+\friction^2)\E\left[\|\b{X}_{0}\|^2\right] + d\friction\right)h^3+ \friction(1+\friction^2) h^2\\
            &\qquad+ \int_{0}^{t_{k}+h} \Exp{\|\b{V}_t -V_t\|^2} + \Exp{\|\b{X}_t -X_t\|^2} \dd t
        \end{aligned}
    \end{equation}
    By the derived Lyapunov drift conditions we have that 
    \begin{equation}
        \disca\coloneqq \sup_{\step>0,\, k\in \N}
            \leq \left( (1+\friction^2+\friction^4)\E\left[\|\b{V}_{t_k}\|^2\right] + (2+\friction^2)\E\left[\|\b{X}_{t_k}\|^2\right] + d\friction \right)\frac{h^3}{3} + \friction(1+\friction^2)\frac{h^2}{2} <\infty
    \end{equation}
    Using $\Exp{\|Z_{0} - \b{Z}_{0}\|^2} = 0$ and a telescope sum we obtain
    \begin{equation}
        \begin{aligned}
            \Exp{\|Z_{t_{k}+h} - \b{Z}_{t_{k}+h}\|^2}
            &=\Exp{\|Z_{t_{k}+h} - \b{Z}_{t_{k}+h}\|^2} - \Exp{\|Z_{t_k} - \b{Z}_{t_k}\|^2} \\
            &\qquad + \sum_{i=0}^{k-1}\Exp{\|Z_{t_{i+1}} - \b{Z}_{t_{i+1}}\|^2} - \Exp{\|Z_{t_i} - \b{Z}_{t_i}\|^2} \\
            &\leq (k+1)\disca \step^2 + \discb \int_{0}^{t_{k}+h} \Exp{\|Z_{t} - \b{Z}_{t}\|^2} \dd t.
        \end{aligned}
    \end{equation}
    Then Grönwall's lemma yields for any $k\in\N$
    \begin{equation}
        \Wc_2^2(\pi R_\step^k,\pi P_{k\step})\leq\Exp{\|Z_{t_{k}} - \b{Z}_{t_{k}}\|^2}\leq k \disca \step^2\expo{\discb k\step}
    \end{equation}
    concluding the proof.
\end{proof}

\section{Relation to Over-relaxed Gibbs Sampling}\label{ssec:gibbs}
In this section, we want to point out the remarkable relation between~\eqref{M:eq:ila} and the (over-relaxed) Gibbs sampler\footnote{
    Indeed, contrary to the presentation in the paper, this connection was the initial inspiration for~\eqref{M:eq:ila}.
}~\cite{adlerOverrelaxationMethodMonte1981,foxAcceleratedGibbsSampling2017}.
Assume for now that our aim is to draw samples from $Z\sim p(z) = \Nc(z \mid 0,\Sigma)$ with some symmetric, positive-definite matrix $\Sigma$.
We can, without loss of generality, assume the mean of the distribution to be zero\footnote{
    In fact, samples from a non-centered Gaussian can be obtained by sampling from the centered one and afterwards adding the mean.
}.
Further, let us separate $Z$ into two blocks $Z=(X,Y)$ with $X,Y\in\R^d$.
A two-block Gibbs sampler for drawing from $Z$ is denoted in \cref{algo:gibbs}~\cite[Section 9]{robert1999monte}.
\begin{algorithm}
    \begin{algorithmic}[1]
        \For{$k=1,2,\dots$}
            \State $X^{k+1} \sim p_{X|Y}(\cdot\mid Y^k)$\label{line:gibbs1}
            \State $Y^{k+1} \sim p_{Y|X}(\cdot\mid X^{k+1})$\label{line:gibbs2}
        \EndFor
    \end{algorithmic}
    \caption{Two-block Gibbs sampler}
    \label{algo:gibbs}
\end{algorithm}
Since $p(z)$ is normal, its conditionals are normal as well~\cite[Section 3.2]{bishopDeepLearningFoundations2024}.
In particular, the updates in lines \ref{line:gibbs1} and \ref{line:gibbs2} of \cref{algo:gibbs} take on the form
\begin{equation}\label{eq:gibbs1}
    \begin{cases}
        X^{k+1} \sim \Nc\bigl( \E[X\mid Y=Y^k], \Sigma_{X|Y}\bigr),\\
        Y^{k+1} \sim \Nc\bigl(\E[Y\mid X=X^{k+1}], \Sigma_{Y|X}\bigr),
    \end{cases}
\end{equation}
where $\Sigma_{X|Y}, \Sigma_{Y|X}$ denote the respective conditional covariance matrices which are, indeed, constant across iterations~\cite[Section 3.2.4]{bishopDeepLearningFoundations2024}.
As shown in \cite{foxAcceleratedGibbsSampling2017}, Gibbs sampling is deeply connected to matrix splitting methods for solving linear equations.
In particular, basic Gibbs sampling corresponds to a stochastic version of the Gauss-Seidel method~\cite{seidel1874uber}.
It is, therefore, a natural question whether acceleration techniques for solving linear systems can be adapted to Gibbs sampling.
This question has already been answered positively in \cite{adlerOverrelaxationMethodMonte1981,foxAcceleratedGibbsSampling2017,nealSuppressingRandomWalks1995}.
Specifically, we may apply \gls{sor} to the Gibbs sampler, leading to the updates
\begin{equation}\label{eq:gibbs2}
    \begin{cases}
        X^{k+1} \sim \Nc\bigl( (1-\omega)X^k + \omega \E[X\mid Y=Y^k], \omega(2-\omega)\Sigma_{X|Y}\bigr),\\
        Y^{k+1} \sim \Nc\bigl( (1-\omega)Y^k + \omega \E[Y\mid X=X^{k+1}], \omega(2-\omega)\Sigma_{Y|X}\bigr),
    \end{cases}
\end{equation}
where the parameter $\omega > 0$ controls the degree of over-relaxation, with larger values typically associated with accelerated convergence.

While the elaboration to this point was restricted to the Gaussian setting, we would like to adapt this acceleration technique to \gls{ula}.
Recall that the standard \gls{ula} update reads as
\begin{equation}
    X^{k+1} = X^k - \tau \nabla U(X^k) + \sqrt{2\tau}N^k,\quad N^k\sim\Nc(0,\id).
\end{equation}
In order to obtain an algorithm resembling the structure of the Gibbs sampler, we introduce an auxiliary variable $Y^k$ and consider the split update 
\begin{equation}\label{eq:aux_gibbs}
    \begin{cases}
        Y^{k+1} = X^k - \tau \nabla U(X^k) + \sqrt{2\sigma_1^2}N^k_1,\\
        X^{k+1} = Y^{k+1} + \sqrt{2\sigma_2^2}N^k_2,
    \end{cases}
\end{equation}
where $(N_i^k)_{i,k}$ are \iid~standard normal and we require $\sigma_1^2 + \sigma_2^2 = \tau$.
It is easy to see that the sequence $(X^k)_k$ in~\cref{eq:aux_gibbs} is equal in distribution to the one generated by \gls{ula}.
Note that, in general, there may not exist a joint density $p(x,y)$ for which \cref{eq:aux_gibbs} constitutes its Gibbs sampler (\cf~\cite{arnoldCompatibleConditionalDistributions1989}).
Nonetheless, we ignore this for now and apply the \gls{sor} over-relaxation to \cref{eq:aux_gibbs}, leading to the scheme
\begin{equation}\label{eq:ula-split-sor}
    \begin{cases}
        Y^{k+1} = (1-\omega)Y^k + \omega (X^k - \tau \grad U(X^k)) + \sqrt{\omega(2-\omega)2\sigma_1^2}N^k_1,\\
        X^{k+1} = (1-\omega)X^k + \omega Y^{k+1} + \sqrt{\omega(2-\omega)2\sigma_2^2}N^k_2.
    \end{cases}
\end{equation}
Let us now reduce this scheme to an iteration solely on $X^k$ again.
Inserting the update rule of $Y^k$ into the second line
and expressing $Y^{k} = \omega^{-1}\bigl(X^{k} - (1-\omega)X^{k-1} - \sqrt{2 \omega (2 - \omega)\sigma_2^2} N_2^{k-1}\bigr)$ in terms of the iterates $X^k$, $X^{k-1}$ and the old noise $N_2^{k-1}$, yields
\begin{equation}
    \begin{aligned}
        X^{k+1} &= (1-\omega)X^k +  (1 - \omega) \Bigl(X^{k} - (1-\omega)X^{k-1} - \sqrt{ 2\omega (2 - \omega)\sigma_2^2} N_2^{k-1}\Bigr) \\
        &\quad+ \omega^2 \bigl(X^k - \tau \grad U (X^k) \bigr) + \omega\sqrt{2 \omega (2 - \omega)\sigma_1^2} N_1^k + \sqrt{2 \omega (2 - \omega)\sigma_2^2} N_2^k.
    \end{aligned}
\end{equation}
Now, choosing $\sigma_2=0$ and $\sigma_1^2 = \tau$ and rearranging terms leads us to
\begin{equation}
    \begin{aligned}
        X^{k+1} = X^k + (1-\omega)^2(X^k-X^{k-1}) - \omega^2\tau \grad U (X^k) + \omega\sqrt{ 2\omega (2 - \omega)\tau} N_1^k.
    \end{aligned}
\end{equation}
Introducing a new parameter $\beta=(1-\omega)^2$ and a new step size variable $\gamma=\omega^2\tau$ and noting that  $\omega(2-\omega) = 1-(1-2\omega + \omega^2) = 1-\beta$ we arrive at 
\begin{equation}\label{eq:gibbs-ila}
    \begin{aligned}
        X^{k+1} = X^k + \beta(X^k-X^{k-1}) - \gamma \grad U (X^k) + \sqrt{2\gamma (1-\beta)} N_1^k
    \end{aligned}
\end{equation}
which is precisely our sampling scheme~\eqref{M:eq:ila}.
While we want to emphasize once more that this derivation is not theoretically well-grounded, as \cref{eq:aux_gibbs} is not the Gibbs sampler of a normal distribution as required in~\cite{foxAcceleratedGibbsSampling2017}, we nonetheless point out that, perhaps surprisingly, the derivation leads precisely to the required trade-off between inertia and randomness rendering \cref{eq:gibbs-ila} a discretization of kinetic Langevin dynamics.

\section{TV Denoising}\label{sec:tv-denoising}
\begin{figure}
    \centering
    \includegraphics[]{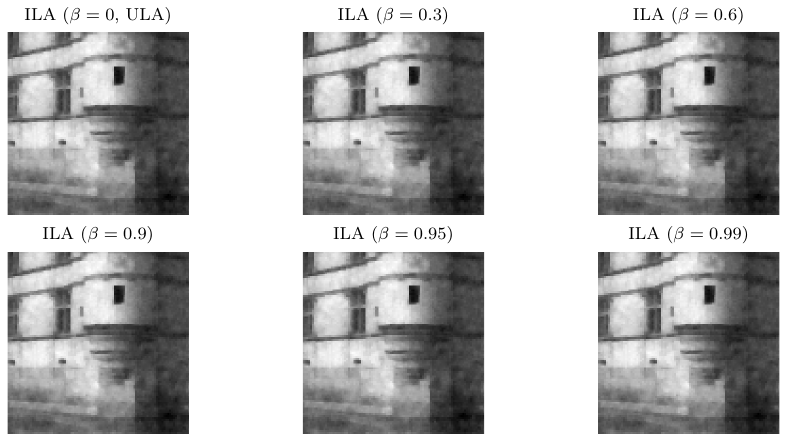}
    \caption{The produced \gls{mmse} estimates.}
    \label{fig:ex-tv-denoising-mmses}
\end{figure}

The \gls{mmse} estimates obtained from the considered \gls{ila} instances can be seen in~\cref{fig:ex-tv-denoising-mmses}.

\section{Training Details for Molecular Structure Generation}\label{app:ebm_training}

All models are trained with contrastive divergence using a replay buffer of size $\num{320}$.
For each iteration, we generate $\num{16}$ new samples initialized with atom positions from a zero-centered normal distribution with a variance matching the data distribution.
To generate samples from the data distribution, we propagate the random samples for $K$ Langevin steps.
We complement the generated samples by adding $\num{48}$ samples from the replay buffer, forming a batch of $\num{64}$ negative samples per iteration.
To account for past changes in the energy surface, we propagate the buffer samples for $\num{15}$ iterations before calculating the contrastive loss.
We apply soft $\ell_2$-regularization on both positive and negative energies to prevent the energy values from drifting to numerically unstable values.
The optimization is performed using the AdamW\footnote{See also: \url{https://docs.pytorch.org/docs/stable/generated/torch.optim.AdamW.html\#adamw}} optimizer with a weight decay factor $\num{1e-2}$.
We conduct a hyperparameter search at $K=200$ Langevin steps over $7$ different step sizes in the interval $[\num{1e-4}, \num{5e-3}]$ for both samplers.
For \gls{ila}, we additionally sweep over three different momentum parameters per step size, chosen according to the theoretical stability bounds derived under the assumption of a Lipschitz constant of $L=1$. We identify the best-performing step size for \gls{ula} with $\num{7.5e-4}$ and for \gls{ila} with $\num{1e-4}$ and a momentum parameter of $\beta=0.987$. This parameter setting is employed throughout \cref{M:ssec:ex-mol-gen}.

\bibliographystyle{siamplain}
\bibliography{references}

\makeatletter\@input{xx.tex}\makeatother